\documentclass{article}

\usepackage{geometry}
\usepackage{graphicx}
\usepackage[dvipsnames]{xcolor}
\usepackage{mathrsfs}
\usepackage{amsfonts}
\usepackage{amssymb}
\usepackage{amsmath}
\usepackage{dsfont}
\usepackage{rotating}
\usepackage{multirow}
\usepackage{hyperref}
\usepackage{subfigure}
\usepackage{comment}
\usepackage{cancel}
\usepackage[noabbrev,capitalise]{cleveref}

\DeclareMathAlphabet\mathbfcal{OMS}{cmsy}{b}{n}  

\newcommand{\x}{\mathbf{x}}
\newcommand{\q}{\mathbf{q}}

\newcommand{\dt}{\Delta t}
\newcommand{\R}{\mathds{R}}

\newcommand{\U}{\mathbf{U}} 
\newcommand{\F}{\mathbf{F}}

\newcommand{\n}{\mathbf{n}}

\renewcommand{\u}{\mathbf{u}}
\renewcommand{\v}{\mathbf{v}}

\renewcommand{\vec}[1]{\mathbf{#1}}

\renewcommand{\div}{\operatorname{div}}
\newcommand{\G}{\mathbf{G}}
\newcommand{\scp}[2]{\left\langle{#1,\, #2}\right\rangle}
\newcommand{\xbi}{{\x}_{\mathbf{b}_i}}  

\newcommand{\inner}[1]{\left\langle #1 \right\rangle}
\newcommand{\urp}{\u^{n+1}_{h,{\gamma^n}}}
\newcommand{\un}{\u^{n}_{h}}
\newcommand{\ud}{\Delta \u_{h}}

\newcommand{\MR}[1]{\textcolor{magenta}{ #1}}

\newtheorem{theorem}{Theorem}[section]
\newtheorem{remark}[theorem]{Remark}

\newtheorem{example}[theorem]{Example}

\newcommand{\diff}[1]{{\mathrm{d}{#1}}}

\title{High order entropy preserving ADER scheme}
\author{Elena Gaburro\thanks{Corresponding author. Inria, Univ. Bordeaux, CNRS, Bordeaux INP, IMB, UMR 5251, 200 Avenue de la Vieille Tour, 33405 Talence cedex, France (\href{mailto:elena.gaburro@inria.fr}{elena.gaburro@inria.fr})}\,, Philipp \"Offner\thanks{Institut f\"ur Mathematik, Johannes Gutenberg Universit\"at, Staudingerweg 9, 55099 Mainz, Germany (\href{mailto:poeffner@uni-mainz.de}{poeffner@uni-mainz.de})}\,, Mario Ricchiuto\thanks{Inria, Univ. Bordeaux, CNRS, Bordeaux INP, IMB, UMR 5251, 200 Avenue de la Vieille Tour, 33405 Talence cedex, France (\href{mailto:mario.ricchiuto@inria.fr}{mario.ricchiuto@inria.fr})}\, and Davide Torlo\thanks{SISSA Mathlab, Mathematics Area, SISSA, via Bonomea 265, 34136, Trieste, Italy (\href{mailto:davide.torlo@sissa.it}{davide.torlo@sissa.it})}}

\begin{document}
	
	\maketitle
	
	%
\begin{abstract}
	In this paper, we develop  a fully discrete entropy preserving ADER-Discontinuous Galerkin (ADER-DG) method. To obtain this desired result, we equip the space part of the method with entropy correction terms that balance the entropy production in space, inspired by the work of Abgrall. Whereas for the time-discretization we apply the relaxation approach introduced by Ketcheson that allows to modify the timestep to preserve the entropy to machine precision. Up to our knowledge, it is the first time that a provable fully discrete entropy preserving ADER-DG scheme is constructed. We verify our theoretical results with various numerical simulations.   
\end{abstract}	
	%
Keywords: Hyperbolic Conservation Laws, ADER, Discontinuous Galerkin, Entropy Conservation/Stability, Structure Preserving, Relaxation Method


\section{Introduction} \label{sec.intro}


In last years, the development of arbitrary high-order entropy preserving discretizations
has attracted a lot of attention in community developing numerical methods for 
hyperbolic PDEs. 
Following the seminal work by Tadmor \cite{tadmor1987numerical,tadmor_2003}, several ideas have surfaced providing general strategies to have control
on the discrete production of entropy. Many works were focused on finding an entropy stable semidiscretization (in space), isolating the entropy conservative flux. These techniques have been applied to many applications, \textit{inter alia} for shallow water equations \cite{gassner2016well,wintermeyer2017entropy,parisot2019entropy}, for Euler's equation \cite{puppo2011numerical,gassner2016split,chandrashekar2013kinetic,ray2016entropy, chandrashekar2013entropy, crean2018entropy}, Navier Stokes problems \cite{gassner2018br1,yamaleev, fernandez2020entropy, manzanero2020entropy,manzanero2020entropy2}, magneto hydro-dynamics problems \cite{chandrashekar2016entropy}, multiphase or multicomponent problems \cite{coquel2021entropy,renac2021entropy,marmignon2022energy} and generally for hyperbolic problems \cite{abgrall2018general, abgrall2020analysis, kuzmin2021entropy, kuzmin2020algebraic, chen2017entropy, carpenter2014entropy, fisher2013high, fisher2013discretely, kopriva2022theoretical, hillebrand2022comparison, renac2019entropy}, also in the Lagrangian framework \cite{carre2009cell,duan2021entropy,chan2021positivity}.
More recently also fully discrete entropy stable or entropy conservative schemes have been introduced, e.g. \cite{ketcheson2019relaxationRK, ranocha2020relaxationRKEuler, abgrall2021relaxation, abgrall2022reinterpretation, friedrich2019entropy, kuzmin2022limiter}. These techniques exploit quite different underlying mechanisms, e.g.  artificial viscosity, limiting techniques, the summation-by-parts (SBP) framework, entropy correction terms, relaxation ansatzes or combination of those. 
Even if the entropy dissipative methods  differ significantly in their underlying ideas, they have demonstrated excellent numerical properties in numerical simulations. \\
The ADER approach has proved to be an accurate numerical scheme in numerous engineering applications. It has also been shown to be
an appropriate setting to design discretizations 
able to preserve different physical properties. First proposed in \cite{toro2002ader, toro2001towards} as an high-order extension of the classical Godunov scheme, the method has been dramatically extended and re-interpreted in different ways, among which we recall
\cite{Dumbser20088209, PCSPH, euromech03, boscheri2014high, ALEDG, boscheri2022cell, ADERGRMHD, dumbser2020glm, dematte2020ader, gaburro2021posteriori, chiocchetti2021high}. However, up to this point no provable high-order entropy conservative/dissipative ADER method can be found in the literature. 
Entropy conservation/stability   is  however essential to
ensure that our numerical approximations converge to the physical relevant weak entropy solution.  In the following work, we will consider this fact and develop for the first time an entropy preserving ADER-DG method. Therefore, we adapt two general techniques to the ADER-DG framework. 
Inspired by \cite{abgrall2018general, abgrall2022reinterpretation}, we handle the space-discretization of our ADER-DG method  via entropy correction terms whereas the entropy production in time is balanced through the relaxation approach proposed in \cite{ketcheson2019relaxationRK,ranocha2020relaxationMultistep,ranocha2020relaxationRKEuler}.\\
The 
paper is organized as follows. In Section~\ref{sec.hyperbolic}, we briefly 
describe hyperbolic conservation laws and the concept of entropy. Then, in Section~\ref{sec.DG}, the ADER-DG method is introduced in its nowadays used explicit space-time Finite Element (FE) interpretation. In Section~\ref{sec.relaxation_DG}, we describe the applied entropy correction techniques and the relaxation approach used in our ADER-DG framework resulting in a provably entropy preserving method. Next, in Section~\ref{sec.test} we verify our theoretical results via various numerical simulations. A conclusion in Section~\ref{sec.concl} finishes this paper.

\section{Hyperbolic Conservation Laws} \label{sec.hyperbolic}


Various natural phenomena and engineering problems can be described by the following system of hyperbolic conservation laws
\begin{equation}
	\label{eq.generalform}
	\partial_t\u+ \operatorname{div} \F (\u)= 0,
\end{equation}
where $\u: \Omega \times \R^+ \to \Psi \subset \R^m$ denotes the conserved variable depending on the space coordinate $\x \in \Omega \subset \R^D $ and time coordinate $t \in R^+$, $\F:\Psi\to \R^{D\, \times\, m}$ is the flux function and $\operatorname{div}$ is the divergence operator. As usual, system  \eqref{eq.generalform} is  equipped with appropriate boundary and initial conditions.  
In the scalar case we 
will use the symbol $u$ for the conservative variable. 
The paper focuses on the two-space dimensional problems \eqref{eq.generalform}, i.e. $\Omega \subset \R^2$. 
However note that all the ideas presented  readily apply to the general case.  In the simulation section, we will show test on scalar linear advection equation, shallow water equations and gas dynamic Euler's equations. All the details on the equations are described in that section.\\
%
%
Due to the non-uniqueness of weak solutions of \eqref{eq.generalform}, the concept of entropy as an additional admissible criterion is necessary.
Following the classical paper \cite{harten1983symmetric}, 
$\eta \colon  \Psi \to \R$ is a convex entropy function with related entropy flux $\G$ and entropy variable $\v= \partial_\u \eta$ verifying $\scp{\v}{\partial_\u  \F_d }=\partial_\u \mathbf{G}_d$.
Smooth solutions of \eqref{eq.generalform}  fulfill additionally  the entropy conservation law
\begin{equation}
	\label{eq.generalform_entropy_conservation}
	\partial_t \eta + \div  \G = 0,
\end{equation}
while admissible weak solutions are  characterized by the entropy inequality 
\begin{equation}
	\label{eq.generalform_entropy_inequality}
	\partial_t \eta + \div \G \leq 0.
\end{equation}
If $\eta$ is  strictly convex, the mapping between the entropy variables $\v =\partial_\u \eta$ and the conservative variables $\u$ is one-to-one. We can either work directly with $\v$ or  $\u$ when solving \eqref{eq.generalform}. Moreover, it holds that the Hessian $\partial_{\u\u}\eta$ (or Jacobian $\partial_\u \v$) is a symmetric positive definite (SPD) matrix and we denote
its inverse by $A_0:= \partial_{\v}\u.$
Finally, our numerical approximations will be constructed in a way to solve 
\eqref{eq.generalform} and fulfill the entropy conditions \eqref{eq.generalform_entropy_conservation} or \eqref{eq.generalform_entropy_inequality} in a global manner. Introducing the total entropy $\mathcal{E}$ defined by the integral
\begin{equation}
	\label{eq:total_U}
	\mathcal{E}(\u)=\int\limits_{\Omega}\eta(\u)\diff \x ,
\end{equation}
total entropy conservation stems immediately from \eqref{eq.generalform_entropy_conservation} and reads
\begin{equation}
	\label{eq.global_entropy_conservation}
	\partial_t \mathcal{E}(\u) + \int\limits_{\partial\Omega}  \G\cdot\hat{\mathbf{n}} \;\diff S = 0.
\end{equation}
Our objective is to provide a strategy allowing to design ADER-DG schemes capable to mimic the last identity  within machine precision.

\section{High order ADER Discontinuous Galerkin schemes} \label{sec.DG}


\subsection{Notation}

This section introduces the basic notation used throughout the paper. 
As mentioned before, we follow the recent interpretation of ADER-DG as an explicit space-time FE discretization. 
We divide the spatial domain $\Omega$ into cells $\Omega_i$ (in our case triangles).
The numerical solution is represented by piecewise polynomials of degree $N \geq 0$ at time level $t^n$ and it is defined by 
\begin{equation}
	\mathbf{u}_h^n(\x,t^n) = \sum \limits_{\ell=0}^{\mathcal{N}-1} \phi_\ell(\x) \, \hat{\mathbf{u}}^{n}_{\ell,i} 
	, \qquad \x \in \Omega_i,
	\label{eqn.uh}
\end{equation}
where $ \phi_\ell(\x) $ is a modal spatial basis and 
$\mathcal{N}$ denotes the total number of expansion coefficients  $\hat{\mathbf{u}}^{n}_{l}$.
The basis function $\phi_\ell(\x) $ is defined through a Taylor series of degree $N$ in the cell $\Omega_i$. Let $\xbi$ be the barycenter of $\Omega_i$ and $h_i$ be the radius of the circumscribed circle around the triangle $\Omega_i$. The basis in $\Omega_i$ is given by  
\begin{equation} 
	\label{eq.Dubiner_phi_spatial}
	\phi_\ell(\x) = \frac{(x - x_{\mathbf{b}_i})^{p_\ell}}{p_\ell! \, h_i^{p_\ell}} \, \frac{(y - y_{\mathbf{b}_i})^{q_\ell}}{q_\ell! \, h_i^{q_\ell}}, \quad 
	\ell = 0, \dots, \mathcal{N}-1, \quad  0 \leq p_\ell + q_\ell \leq N \quad \text{with} \quad \mathcal{N}= \frac{1}{d!} \prod \limits_{m=1}^{d} (N+m).
\end{equation} 
\begin{remark}[Extension to the general $P_NP_M$ framework]
	If we choose $N=0$, we obtain as usual a first order FV method. Instead of increasing now the order of accuracy in one cell by setting $N>0$, one can,  for instance, apply a higher-order FV reconstruction up to order $M+1$. 
	Therefore, our method is included in the more general class of $P_NP_M$ schemes proposed 
	in~\cite{Dumbser20088209} and further developed in \cite{gaburro2021posteriori,gaburro2020high} which unifies higher-order FV and DG schemes into a common framework. Since the application of the entropy correction term inside a  FV scheme is more problematic as described in \cite{abgrall2022reinterpretation}, we restrict ourself in this paper on the ADER-DG approach only. Extensions to general $P_NP_M$ are planned for future  work. 
\end{remark}

Differently from classical Runge-Kutta discontinuous Galerkin method, the ADER-DG approach is not build on a method of line (MOL) ansatz but on a local predictor-corrector approach. The advantage of the predictor-corrector procedure is that the information exchange between the elements is done only once at each time step, only during the corrector step. This results in advantages in runtime and efficiency of the scheme \cite{balsara2013efficient}.
Therefore, we also define a space-time basis which  locally approximates our predictor $\q_h$ inside the element $C_i^n=\Omega_i \times [t^n,t^{n+1}]$. It is given by 
\begin{equation} 
	\q_h^n(\x, t) = \sum_{\ell=0}^{\mathcal{Q}-1} \theta_\ell (\x, t) \hat{\q}_\ell^n, \quad (\x,t) \in C_i^n,  
	\label{eqn.qh}
\end{equation} 
where  
\begin{equation}
	\begin{aligned}
		\label{eq.Dubiner_phi}
		& \theta_\ell(x,y,t)|_{C_i^n} = \frac{(x -  x_{\mathbf{b}_i})^{p_\ell}}{{p_\ell}! \, h_i^{p_\ell}} \, \frac{(y -  y_{\mathbf{b}_i})^{q_\ell}}{{q_\ell}! \, h_i^{q_\ell}}
		\, \frac{(t - t^n)^{r_\ell}}{{r_\ell}! \, h_i^{r_\ell}}, \\
		& \ell = 0, \dots, \mathcal{Q}, 
		\quad 0 \leq p_\ell + q_\ell + r_\ell \leq N \quad \text{with} \quad \mathcal{Q}= \frac{1}{(d+1)!} \prod \limits_{m=1}^{d+1} (N+m)
	\end{aligned}
\end{equation} 
is our \textit{modal space--time} basis for the polynomials of degree $N$. Here, $\mathcal{Q}$ is the number of degrees of freedom in one space-time cell $C_i$.

\subsection{ADER-DG method}\label{subsec.ADER-DG}

ADER-DG schemes consist in two steps for the update of one timestep: the local predictor and the corrector, which we describe now separately.
\subsubsection{Local Predictor Step}
The method starts by calculating a high-order local predictor solution of the underlying PDE \eqref{eq.generalform}  using a classical Galerkin approach inside each element $C_i$. We obtain  
\begin{equation}
	\int_{t^n}^{t^{n+1}} \!\! \int_{\Omega_i^{\circ} } \theta_k \, \partial_t \q_h \,\diff \x \, \diff t 
	+ \int_{t^n}^{t^{n+1}} \!\! \int_{\Omega_i^{\circ} } \theta_k \,\div \F(\q_h) \, \diff \x \,\diff t
	= \mathbf{0}
	\label{eq.predictorEq}
\end{equation}
with  $\Omega_i^{\circ} = \Omega_i \backslash \partial \Omega_i$ being the interior of the spatial cell.\\
Applying integration by parts in time yields finally 
\begin{equation}
	\begin{aligned}
		& \int_{\Omega_i^{\circ}}  \theta_k(\x,t^{n+1}) \q_h(\x,t^{n+1}) \, \diff \x -
		\int_{\Omega_i^{\circ}} \theta_k(\x,t^{n}) \u_h(\x,t^{n}) \, \diff \x - \\
		&\int_{t^n}^{t^{n+1}} \!\!\!  \int_{\Omega_i^{\circ} } \!\!\! \partial_t \theta_k(\x,t) \q_h(\x,t)  \diff \x \,\diff t + 
		\int_{t^n}^{t^{n+1}} \!\!\!  \int_{\Omega_i^{\circ} } \!\!\! \theta_k(\x,t) \div\F(\q_h(\x,t))  \diff \x \,\diff t
		= \mathbf{0},
		\label{eq:DOFpredictor}
	\end{aligned}
\end{equation}
in which now $\u_h(\x,t^{n})$ is given. Equation \eqref{eq:DOFpredictor} involves only the local polynomial $\q_h^n$ within a single  space-time element.
The space-time  predictors \eqref{eq:DOFpredictor} can be thus solved independently in each element. In practice,
we approximate their solution 
by a simple fixed-point iteration procedure (discrete Picard iterations) as follows
\begin{equation}
	\begin{aligned}
		& \int_{\Omega_i^{\circ}}  \theta_k(\x,t^{n+1}) \q_h^{(k+1)}(\x,t^{n+1}) \, \diff \x -
		\int_{t^n}^{t^{n+1}} \!\!\!  \int_{\Omega_i^{\circ} } \!\!\! \partial_t \theta_k(\x,t) \q_h^{(k+1)}(\x,t)  \diff \x \,\diff t =\\
		&\int_{\Omega_i^{\circ}} \theta_k(\x,t^{n}) \u_h(\x,t^{n}) \, \diff \x + 
		\int_{t^n}^{t^{n+1}} \!\!\!  \int_{\Omega_i^{\circ} } \!\!\! \theta_k(\x,t) \div\F(\q_h^{(k)}(\x,t))  \diff \x \,\diff t.
		\label{eq:DOFpredictor1}
	\end{aligned}
\end{equation} 
For further details one can refer to \cite{gaburro2020high,gaburro2021unified} where the same notation has been used.

\subsubsection{Corrector Step}

Once the  local space-time prediction of the  solution  $\q_h^n(\x, t) $ has been computed, we can write a one step  update in time as follows. We
start by multiplying \eqref{eq.generalform} with  a spatial-only test function  $\phi_k$ (cf. equation \eqref{eq.Dubiner_phi_spatial}), and integrate in space and time: 
\begin{equation}
	\int_{t^n}^{t^{n+1}} \int_{\Omega_i} \phi_k \left(\partial_t \u  + \div \F(\u) 
	\right)   \diff \x \,\diff t
	= \mathbf{0}.
	\label{eqn.pde.weak}
\end{equation}
Inserting now \eqref{eqn.uh} in the time evolution variable and our local information \eqref{eq:DOFpredictor} in the flux, we can apply the divergence theorem in \eqref{eqn.pde.weak}, and we get our final update step of
the ADER-DG method:
\begin{equation}
	\begin{aligned}
		& \left( \, \int_{\Omega_i} \phi_k \phi_l \, d\x\right)
		\left( \hat{\u}_\ell^{n+1} - \hat{\u}_\ell^{n} \, \right) 
		+ \int_{t^n}^{t^{n+1}} \!\! \int_{\partial \Omega_i} \!\! \phi_k \mathcal{F}\left(\q_h^-, \q_h^+ \right) \cdot \mathbf{n} \, \diff S \, dt \ - 
		\int_{t^n}^{t^{n+1}} \!\! \int_{\Omega_i } \nabla \phi_k \cdot \F(\q_h)  \diff \x \,\diff t 
		\, =\,  \mathbf{0},
		\label{eq:ADER-DG}
	\end{aligned}
\end{equation}
where $\mathcal{F}$ is the numerical flux, and $\mathbf{n}$ the external normal vector at the boundaries of the cell. Since we have the local predictor information $\q_h$, we can apply high-order quadrature rules for the approximation of time-integrals in \eqref{eq:ADER-DG} and we finally obtain our update for the solution at $t^{n+1}$. We remark that the information exchange in-between the elements plays a role only in the corrector step. 
As mentioned before, the presented interpretation of ADER-DG is an explicit space-time FE method.

\section{Fully discrete entropy preserving ADER discontinuous Galerkin schemes} \label{sec.relaxation_DG}


As described in the previous section, the ADER-DG approach consists of two steps:
a set of  local predictors, and a fully explicit global corrector. The corrector step enforces
conservation, and thus consistency with weak solutions of \eqref{eq.generalform}.
However, the entropy (in)equality \eqref{eq.generalform_entropy_conservation} (or \eqref{eq.generalform_entropy_inequality}) might not be guaranteed. 
Due to the  inherent space-time nature of the method, to enforce this additional constraint 
we need to have simultaneous control  of the  entropy production in space and in time.  
To this end we extend  two 
techniques  previously proposed in literature: 
\begin{enumerate}
	\item for the balance of the spatial discretization in our ADER-DG method \eqref{eq:ADER-DG}, we apply entropy correction terms first proposed by Abgrall in \cite{abgrall2018general},
	\item for handling the time-integration we deal with the relaxation approach described in \cite{ketcheson2019relaxationRK, ranocha2020relaxationMultistep}.
\end{enumerate}
In the following section 
we explain how these techniques can be adapted to \eqref{eq:ADER-DG}. 

\subsection{Entropy correction in space}

The approach proposed by Abgrall in \cite{abgrall2018general} is based on the following idea. We introduce
a simple entropy correction term to the underlying discretization which is lead by a parameter which varies in every cell. 
The parameter is chosen in order to locally balance the entropy production at every degree of freedom and, simultaneously, not to violate the conservation property of the scheme. 
It can be interpreted as the solution of an optimization problem \cite{abgrall2022reinterpretation} or as an additional viscosity/anti-viscosity term. 
Inspired by \cite{abgrall2018general}, we introduce an entropy correction term with artificial 
diffusion into the corrector of the ADER formulation \eqref{eq:ADER-DG} using a generic element $\psi$ from our underlying approximation space.  It yields
\begin{equation}\label{eq:ADERupdateDiffusion}
	\begin{aligned}
		& \int_{\Omega_i} \psi  (\x)
		\left( \u_h^{n+1}(\x) - \u_h^{n}(\x) \, \right) \, \diff \x 
		+ \int_{t^n}^{t^{n+1}} \!\! \int_{\partial \Omega_i} \!\! \psi(\x) \mathcal{F}_{\mathbf{n}}\left(\q_h^-(\x,t), \q_h^+(\x,t) \right) \, \diff S \, \diff t \ - \\
		& \int_{t^n}^{t^{n+1}} \!\! \int_{\Omega_i } \!\!\! \nabla \psi(\x)\cdot \F(\q_h(\x,t))  \diff \x \,\diff t 
		+   \int_{t^n}^{t^{n+1}}\!\! \alpha_{\Omega_i}\int_{\Omega_i}\nabla \psi(\x) \, A_0\, \nabla \v_h(\q_h(\x,t))  \, \diff \x \,\diff t 
		\, =\,  \mathbf{0},
	\end{aligned}
\end{equation}
where $\v_h$ is the approximation of the entropy variables $\v\in \R^{m}$ and $A_0\in \R^{m\times m}$ is a positive definite matrix\footnote{To obtain dimensionally correct terms, the inverse of the Hessian of the entropy $\eta_{\U\U}$ evaluated in an average value of the cell will be used.}.
Using some quadrature formula for the time integration, from \eqref{eq:ADERupdateDiffusion} we obtain 
\begin{equation}
	\begin{aligned}
		& \int_{\Omega_i} \psi  (\x)
		\left( \u_h^{n+1}(\x) - \u_h^{n}(\x) \, \right) \, d\x 
		+ \Delta t \sum_{s=0}^S\beta_s \int_{\partial \Omega_i} \!\! \psi(\x) \mathcal{F}_{\mathbf{n}}\left(\q_h^{s,-}(\x), \q_h^{s,+}(\x) \right)  \, \diff S  - \\
		& \Delta t \sum_{s=0}^S \beta_s  \!\!\int_{\Omega_i} \nabla \psi(\x)\cdot \F(\q^s_h(\x)) \,\diff \x
		+ \Delta t \sum_{s=0}^S \beta_s \,\alpha^s_{\Omega_i}  \int_{\Omega_i} \nabla \psi(\x) A_0 \nabla \v_h(\q^s_h(\x)) \diff \x
		\, =\,  \mathbf{0},
	\end{aligned}
\end{equation}
where $\beta_s$ are the time quadrature weights.
It is left to define all $\alpha_{\Omega_i}^s$. We choose them in order to have an entropy dissipative (conservative) semi-discretization update.
This means that at each 
subtimestep $t^{n,s}$, substituting $\psi = \v_h(\q^s_h(x))$ and summing over the variables, we should have
\begin{align}
	&\int_{\partial \Omega_i} \!\!  \inner{\v_h(\q^s_h(\x)), \mathcal{F}_{\mathbf{n}}\left(\q_h^{s,-}(\x), \q_h^{s,+}(\x) \right) } \, \diff S  - 
	\!\!\int_{\Omega_i} \nabla  \v_h(\q^s_h(\x)) : \F(\q^s_h(\x)) \,\diff \x +\\
	&\alpha^s_{\Omega_i} \int_{\Omega_i}  \nabla  \v_h(\q^s_h(\x))^T A_0 \nabla \v_h(\q^s_h(\x)) \diff \x
	\, \geq \,  \int_{\partial \Omega_i} \mathcal{G} (\v_h(\q^{s,-}(\x)),\v_h(\q^{s,+}(\x))) \cdot \mathbf{n}\, \diff S,
\end{align}
where $\mathcal{G}$ is a consistent  numerical entropy flux.
In order to achieve this, we split the numerical flux into a central part and a diffusive one as 
\begin{equation}
	\mathcal{F}(\q^+,\q^-)\cdot \mathbf{n} = \mathcal{F}^c_{\mathbf{n}}(\q^+,\q^-)+\mathcal{F}^d(\q^+,\q^-).
\end{equation}
As an example with a {local Lax-Friedrich (Rusanov)} entropy type flux, we would have 
\begin{equation}\label{eq:numericalFluxParts}
	\mathcal{F}^c_{\mathbf{n}}(\q^+,\q^-):= \frac{\mathbf{F}(\q^+)+\mathbf{F}(\q^-)}{2} \cdot \mathbf{n}, \ \text{ and }\
	\mathcal{F}^d(\q^+,\q^-):= -\mathbf{a}_{LF}(\q^+,\q^-) \frac{\q^+-\q^-}{2}
\end{equation}
with $\mathbf{a}_{LF}$ the spectral radius of the flux Jacobian. \\
For simplicity, we introduce the following definitions of the central spatial contribution, the diffusion flux term and the entropy fix  
\begin{align*}
	&\mathtt{F}_{i}^s:= \int_{\partial \Omega_i} \!\! \left\langle \v_h(\q^s_{h,i}(\x)) , \mathcal{F}^c_{\mathbf{n}}\left(\q_h^{s,-}(\x), \q_h^{s,+}(\x) \right) \right\rangle \, dS  \ - \!\! \int_{\Omega_i^{\circ} } \!\!\! \nabla \v_h(\q_h^s(\x)) : \F(\q_h^s(\x)) \,d\x,\\
	&\mathtt{D}_{i}^s:= \int_{\partial \Omega_i} \!\left \langle\v_h(\q^s_{h,i}(\x)),  \mathcal{F}^d\left(\q_h^{s,-}(\x), \q_h^{s,+}(\x) \right) \right\rangle \, dS , \\
	& \mathtt{E}_{i}^s:= \int_{\Omega_i} \nabla \v_h(\q^s_h(\x))^T\, A_0\, \nabla \v_h(\q_h^s(\x)) \,d\x  \geq 0.
\end{align*}
With these definitions, we can finally determine $\alpha_{\Omega_i}^s$ in such a way to ensure entropy conservation or entropy dissipation.
We start focusing on the entropy dissipative case, therefore
a simple choice would be to determine $\alpha_{\Omega_i}^s$  such that the entropy contribution of the central part vanishes, i.e.,
\begin{equation}
	\int_{\partial \Omega_i}\mathcal{G}(\v_h(\q^{s,-}(\x)),\v_h(\q^{s,+}(\x))) \cdot \mathbf{n}\, \diff S= \mathtt{F}_i^s + \alpha_{\Omega_i}^s \mathtt{E}_i^s, \qquad s=1,\dots, S.
\end{equation}
This leads to the definition of 
\begin{equation}\label{eq:alpha_i-s}
	\alpha_{\Omega_i}^s =\frac{\int_{\partial \Omega_i}\mathcal{G}(\v(\q^{s,-}(\x)),\v(\q^{s,+}(\x))) \cdot \mathbf{n}\, \diff S- \mathtt{F}_i^s}{\mathtt{E}_i^s}.
\end{equation}
Since the diffusive part $\mathtt{D}_{i}^s$ decreases always the entropy, we obtain the following estimation
\begin{equation*}
	\begin{aligned}
		\int_{t^n}^{t^{n+1}} \int_\Omega& \v_h(\x,t)^T\partial_t \u_h(\x,t)  \, \diff \x  \approx - \Delta t \sum_{s=0}^S \beta_s \sum_{i\in \Omega} \left(\mathtt{F}_i^s + \alpha_{\Omega_i}^s \mathtt{E}_i^s + \mathtt{D}^s_{i}\right) \\
		&=- \Delta t \sum_{s=0}^S \beta_s \left[  \sum_{\Omega_i \in \Omega}\int_{\partial \Omega_i}     \mathcal{G}(\v(\q^{s,-}(\x)),\v(\q^{s,+}(\x))) \cdot \mathbf{n}\, \diff S + \sum_{i\in \Omega} \mathtt{D}_i^s \right] \\
		&\leq - \Delta t \sum_{s=0}^S \beta_s \int_{\partial \Omega}     \G(\v(\q^s(\x))) \cdot \mathbf{n}\, \diff S.
	\end{aligned}
\end{equation*}
With this definition, we correct the contribution from the centered term of the numerical flux. In particular, we are able to split the update into an entropy conservative part ($\mathtt{F}_i^s+\alpha_{\Omega_i}^s \mathtt{E}_i^s$) and an entropy dissipative part ($\mathtt{D}_i^s$). 
This splitting will be crucial when introducing the relaxation technique, as it allows to choose which part we need to balance also in time. 
\begin{remark}[Division by zero]
	The division performed in \eqref{eq:alpha_i-s} is numerically dangerous as in flat areas $\normalfont \texttt{E}_i^s$ might be 0. In order not to risk to encounter NaN, we perform that division only if $\normalfont \texttt{E}_i^s \geq \varepsilon^s$, otherwise we set $\alpha_{\Omega_i}^s =0$. This tolerance $\varepsilon^s$ is chosen in the following way
	\begin{equation}
		\varepsilon^s : = \bar{h}^{\,N} \max_{i\in \Omega} \normalfont \texttt{E}_{i}^s,
	\end{equation}
	where $N$ is the polynomial degree and $\bar{h}$ is the average mesh size. This guarantees not to introduce extra terms where the solution is constant or almost constant. 
\end{remark}

\begin{remark}[Correction for first order finite volume schemes]
	The special case of $N=0$, i.e., the finite volume method, does not fit in our framework as the terms $\mathtt{E}$ are all zero.
	In this case, we can resort to the classical characterization by Tadmor \cite{tadmor1987numerical,tadmor_2003}.  
	On unstructured grids we introduce  the jumps  $[\! [ (\cdot) ]\!] = (\cdot)^+ -(\cdot)^- $ on a given mesh face with 
	normal $\n_f$. The generalization of the classical Tadmor condition \cite{ray2016entropy, tadmor_2003}
	\begin{equation}\label{eq:Tadmor_FV}
		[\![ \v_h]\!]^T\mathcal{F}^{co}_{\mathbf{n}_f}  = [\![ \pmb{\psi}(u_h)]\!] \cdot \n_f
	\end{equation}
	defines an entropy conservative flux $\mathcal{F}^{co}$ consistent with the numerical entropy flux
	\begin{equation}\label{eq:Tadmor_FV_EC}
		\mathcal{G}_{\n_f}  :=   \overline{\v}_h^T\mathcal{F}^{co}_{\mathbf{n}_f}  - \overline{\pmb{\psi}}(u_h)\cdot \n_f,
	\end{equation}
	having denoted by $\overline{\cdot}$ the arithmetic average, and by $\pmb{\psi}$ the  entropy potential defined by
	\begin{equation}\label{eq:E-potential}
		\G = \v^T\F -\pmb{\psi}.
	\end{equation}              
	A given numerical flux $\mathcal{F}_{\n_f}$ may not verify \eqref{eq:Tadmor_FV}. However, we can use the  face jumps  as a correction term and define
	\begin{equation}\label{eq:FV-correction}
		\mathcal{F}^{co}_{\mathbf{n}_f}=  \mathcal{F}_{\mathbf{n}_f} -  \alpha_f  A_0 (\v^+ - \v^-) =\mathcal{F}_{\mathbf{n}_f} -  \alpha_f  A_0 [\![\v]\!].
	\end{equation}
	By imposing the entropy conservation condition \eqref{eq:Tadmor_FV}, we compute a limiting value 
	%
	\begin{equation}\label{eq:correction-pnpm-j}
		\alpha_0  := \dfrac{ - [\![ \pmb{\psi}]\!] \cdot\n_f  + [\![\v]\!]^T \mathcal{F}_{\mathbf{n}_f}   }{[\![\v]\!]^TA_0[\![\v]\!]} =
		\dfrac{ [\![ \G]\!] \cdot \n_f - [\![\v_h^T\F ]\!] \cdot\n_f  + [\![\v]\!] ^T\mathcal{F}_{\mathbf{n}_f}   }{[\![\v]\!]^TA_0[\![\v]\!]} .
	\end{equation}
	Entropy conservation is ensured  for $\alpha_f=\alpha_0$,  	while  semi-discrete entropy stability is obtained for $\alpha_f > \alpha_0$.
	Note that this  comparison result is relatively classical, see \cite{tadmor_2003}.
	In this work, we will focus on high order methods, and we will not discuss further this case. The authors are anyway considering this formulation to obtain an \textit{a posteriori} subcell finite volume limiter, in order to deal with solutions with strong shocks. This will be object of future studies, as the coupling of the two techniques is delicate and must be handled carefully. 
\end{remark}  


\subsection{Relaxation applied to ADER-DG}
Would the time integration be exact (vis in the time continuous case), the correction introduced in the previous section would allow to satisfy
\emph{exactly} entropy conservation/stability within each element. But, unfortunately, in general this is not the case and moreover 
the fully discrete nature of the ADER approach requires to handle the temporal discretization simultaneously with the spatial one.
To achieve exact entropy conservation thus we adapt 
%
the relaxation approach proposed in \cite{ketcheson2019relaxationRK}, and applied in different contexts 
\cite{abgrall2021relaxation, kang2021entropy, ranocha2020relaxationRKEuler, ranocha2020relaxationMultistep}.
The key idea of the relaxation ansatz is to project the final numerical approximation $ \u^{n+1}_h$ back to the manifold where the entropy is constant (or dissipative), as in  the underlying continuous setting. 
This is obtained multiplying the time-step $\Delta t$ with a scalar quantity $\gamma^n$ and using this new time-step for the update. 
We describe now how $\gamma^n$ can be chosen in order to preserve (or dissipate) the total entropy \eqref{eq:total_U}.
We set 
\begin{equation}
	\Delta \u_h := \u^{n+1}_h -\u^{n}_h, 
\end{equation} 
where $\u^{n+1}_h$ is obtained solving \eqref{eq:ADERupdateDiffusion}. Now, we have to find 
$\gamma^n$ and $\u^{n+1}_{h,{\gamma^n}} := \u_h^n + \gamma^n \Delta \u_h $ such that the total entropy is controlled. We start by estimating the total entropy \eqref{eq:total_U} as
\begin{align}
	\mathcal{E}(\urp) &= \int_{\Omega} \eta(\urp) d\x \approx \int_{\Omega} \eta(\un) d\x+ \gamma^n \langle \eta_u(\un) ,\ud \rangle \diff \x \\
	&\approx \mathcal{E}(\un) - \gamma^n \Delta t \sum_{s=0}^S \beta_s \sum_{i\in \Omega} \left(\mathtt{F}_i^s + \alpha_{\Omega_i}^s \mathtt{E}_i^s + \mathtt{D}^s_{i} \right) \\
	&\leq \mathcal{E}(\un) -\gamma^n\Delta t \sum_{s=0}^S \beta_s \int_{\partial \Omega } G(\v^s) \cdot \mathbf{n} \diff S \approx \mathcal{E}(\un) -\int_{t^n}^{t^{n+1}_{\gamma^n}}  \int_{\partial \Omega} G(\v^s) \cdot \mathbf{n} \diff S.
\end{align}
This leads us to the following (non-)linear equation: 
\begin{equation}\label{eq:ADERgammaEquation}
	\mathcal{E}(\urp) = \int_{\Omega}  \eta(\un +\gamma^n \ud )d\x \stackrel{!}{=} \int_\Omega \eta(\un)\diff \x - \gamma^n\Delta t \sum_{s=0}^S \beta_s \sum_{i\in \Omega}  \left(  \mathtt{F}_{i}^s+\alpha_{\Omega_i}^s\mathtt{E}_{i}^s +\mathtt{D}_{i}^s \right).
\end{equation}
This is equivalent to finding the zero of 
\begin{equation}\label{eq:relaxation_residual}
	R(\gamma^n):= \int_{\Omega}  \eta(\un +\gamma^n \ud )d\x -\int_\Omega \eta(\un)\diff \x + \gamma^n\Delta t \sum_{s=0}^S \beta_s \sum_{i\in \Omega}  \left(  \mathtt{F}_{i}^s+\alpha_{\Omega_i}^s\mathtt{E}_{i}^s +\mathtt{D}_{i}^s \right).
\end{equation}
In case we want to exactly preserve the entropy in time, we can remove the $\mathtt{D}_i^s$ terms from equation \eqref{eq:ADERgammaEquation}.
From \eqref{eq:ADERgammaEquation}, we can finally calculate $\gamma^n$ at every time-step and we update the new value of the approximated solution and the new time-step via 
\begin{equation}\label{eq:update_relax}
	\begin{cases}
		\tilde{\u}^{n+1}_{h,\gamma^n} := \un + {\gamma^n} \Delta \un,\\
		\tilde{t}^{n+1}_{\gamma^n} := t^n + {\gamma^n} \Delta t.
	\end{cases}	
\end{equation}

\begin{remark}[Division by zero (Relaxation)] As discussed in Remark 4.1, when $\mathtt{E}_i^s< \varepsilon$ we set  $\alpha_{\Omega_i}=0$.
	To ensure control of the entropy production up to machine accuracy this particular case needs to be treated also in the relaxation procedure.
	In our approach, if $\mathtt{E}_i^s< \varepsilon$ we replace  $\mathtt{F}_{i}^s+\alpha_{\Omega_i}^s\mathtt{E}_{i}^s $
	by $\int_{\partial \Omega_i} G(\v^s)\cdot \mathbf{n} \diff S$ in  \eqref{eq:ADERgammaEquation}. This guarantees that we still have 
	\begin{equation}
		\mathcal{E}(\urp) = \mathcal{E}(\u^n_h) - \sum_{s=0}^S \beta^s \left[ \int_{\partial \Omega} G(\v^s) \cdot\mathbf{n} \,\diff S + \sum_{\Omega_i \in \Omega} \mathtt{D}_i^s \right],
	\end{equation}
	which is otherwise guaranteed by the definition of  $\alpha_{\Omega_i}$.
\end{remark}

%

It should be pointed out that in case of quadratic entropies \eqref{eq:ADERgammaEquation}  is explicitly solvable as demonstrated in Example \ref{ex.quadratic} below, while for more nonlinear entropies a nonlinear solver for the scalar equation \eqref{eq:ADERgammaEquation} must be used. In our numerical simulations in Section \ref{sec.test}, we use a simple Newton method to solve \eqref{eq:ADERgammaEquation}. 
\begin{remark}[Newton method for \eqref{eq:relaxation_residual}]
	When we want to find the (non trivial) zero of \eqref{eq:relaxation_residual}, the Newton method or variants can help in finding the solution. For Newton method also the derivative of $R$ is required, but it is easily computable as
	\begin{equation}
		\frac{dR(\gamma^n)}{d\gamma^n} = \int_{\Omega}  \left\langle \v(\un +\gamma^n \ud ),\ud \right\rangle d\x  +\Delta t \sum_{s=0}^S \beta_s \sum_{i\in \Omega}  \left(  \mathtt{F}_{i}^s+\alpha_{\Omega_i}^s\mathtt{E}_{i}^s +\mathtt{D}_{i}^s \right).
	\end{equation}
\end{remark}In that case, to compute the Jacobian, the entropy variable definition is, again, necessary.
Finally, we have to stress out that due to the results of \cite{ranocha2020relaxationMultistep}, we can ensure that
if the underlying time method is of order $p\geq2$,
\begin{itemize}
	\item there exists always a unique $\gamma>0$ that solves \eqref{eq:ADERgammaEquation}. This $\gamma$ satisfies $\gamma=1+ \mathcal{O}(\Delta t^{p-1})$ under a reasonable small time-step $\Delta t$, and 
	\item the relaxation solution $\tilde{\u}^{n+1}_{h,\gamma^n} $ is also of order $p$.  
\end{itemize}

Using now both, the entropy correction terms and the relaxation approach, we obtain an entropy preserving ADER-DG method for which we will verify the properties next in the numerical experiments section, but before we finish this part with the following already mentioned example.

\begin{example}[Quadratic entropy]\label{ex.quadratic}
	In case of a quadratic entropy for scalar equations 
	$$
	\eta(u) :=\frac{ \langle u, u \rangle }{2},
	$$
	we can exploit the expression $\mathcal{E}(\urp)$ as
	\begin{align}
		\mathcal{E}(\urp)& = \int_{\Omega} \eta(\un +\gamma^n \ud ) d\x = \frac{1}{2} \int_{\Omega} \langle  \un +\gamma^n \ud ,  \un +\gamma^n \ud \rangle  \\
		&= \frac{1}{2}\int_{\Omega} \langle  \un ,  \un \rangle  + 2 \gamma^n \langle  \ud ,  \un \rangle  + (\gamma^n)^2 \langle \ud ,\ud \rangle d \x  \\
		&= \mathcal{E}(\un) + \int_{\Omega}  \gamma^n \langle  \ud ,  \un \rangle  + \frac{1}{2}(\gamma^n)^2 \langle \ud ,\ud \rangle d \x,
	\end{align}
	so that \eqref{eq:ADERgammaEquation} reduces to
	\begin{align}
		\int_{\Omega}  \gamma^n \langle  \ud ,  \un \rangle  + \frac{1}{2} (\gamma^n)^2 \langle \ud ,\ud \rangle \diff \x \stackrel{!}{=} -\gamma^n\Delta t \sum_{s=0}^S \beta_s \sum_{i\in \Omega}  \left(  \mathtt{F}_{i}^s+\alpha_{\Omega_i}^s\mathtt{E}_{i}^s +\mathtt{D}_{i}^s \right).
	\end{align}
	The solution $\gamma^n=0$ is not the one we search, the other one is
	\begin{equation}\label{eq_gamme_energy}
		\gamma^n =-2 \frac{\Delta t \sum_{s=0}^S \beta_s \sum_{i\in \Omega}  \left(  \mathtt{F}_{i}^s+\alpha_{\Omega_i}^s\mathtt{E}_{i}^s +\mathtt{D}_{i}^s \right)+\int_{\Omega} \langle  \ud ,  \un \rangle \diff \x }{\int_{\Omega} \langle \ud ,\ud \rangle d \x},
	\end{equation}
	which is close to one. 
\end{example}


\section{Numerical results} \label{sec.test}

In this section we test the proposed algorithm against many physical problems: linear advection, shallow water and Euler equations.
We will focus on different aspects in the various tests.
In most cases we will compare the ADER-DG method with the proposed ADER-DG entropy fix with the relaxation technique.
We will study the convergence analysis, the entropy error in time and the value of $\gamma^n$ in time.
For each problem we will introduce the equations, the entropy variables, the entropy fluxes and the matrix $A_0$.

\subsection{Linear advection and rotation}
The first equation we treat is the linear advection equations in two forms.
The equation we consider is 
	\begin{equation}
		\partial_t u + \vec a \cdot \nabla u  = 0,
	\end{equation}
	with $\vec a=(a_1, a_2)$ which can be constant or a rotation vector, i.e., $\vec a = (-y,x)$. The entropy, entropy variable and entropy flux are 
	\begin{equation}
		\eta(u)= u^2/2, \quad v(u) = u,\quad G(u) = \eta(u)\,\vec a.
	\end{equation}
	In this case we have trivially that $A_0 = 1$. 

\subsubsection{Test: traveling bump}\label{sec:test_traveling_bump}
The first smooth test that we consider is given by the following initial condition 
\begin{equation}
	u_0(\vec x) = 
	\begin{cases}
		e^{1-\frac{1}{1-r^2}} & \text{if } r<1,\\
		0&\text{else},	
	\end{cases}
\end{equation}
where $r=\sqrt{x^2+y^2}$, the advection coefficient is $\vec a = (1,0)^T$, the domain is $\Omega=[-1.5,1.5]\times [-1.5,1.5]$ with periodic boundary conditions in the $x$ direction and wall boundary conditions in the $y$ direction. 
The exact solution is given by $u_{ex}(x,t)= u_0(\vec x-\vec a t)$.

\begin{table}
	\caption{Order of convergence for the ADER-DG schemes equipped with the relaxation technique at time $t_f=3.0$ for the $L_2$ norm error of the travelling bump.} 
	\label{tab.TravellingBump}
	\begin{center} 	
		\begin{tabular}{|ccc|ccc|ccc|} 
			\hline
			\multicolumn{3}{|c|}{$P_1 +$ Rel}        &   \multicolumn{3}{|c|}{$P_2 +$ Rel}      &    \multicolumn{3}{|c|}{$P_3 +$ Rel}   \\
			$\Delta x$   &   $L_2(h-h_{ex})$   & $\mathcal{O}(L_2)$ & $\Delta x$   &   $L_2(h-h_{ex})$   & $\mathcal{O}(L_2)$ & $\Delta x$   &   $L_2(h-h_{ex})$   & $\mathcal{O}(L_2)$   \\
			\hline
			2.37E-2 & 3.51E-3 &   -  & 3.31E-2 & 4.47E-4 &  -    & 4.30E-2 & 2.07E-4 &   -    \\
			1.69E-2 & 1.93E-3 & 1.78 & 2.43E-2 & 2.74E-4 & 2.77  & 3.12E-2 & 6.77E-5 &  3.47  \\ 
			1.23E-2 & 1.09E-3 & 1.82 & 1.75E-2 & 1.05E-4 & 2.91  & 2.27E-2 & 2.02E-5 &  3.79  \\
			8.91E-3 & 5.96E-4 & 1.84 & 1.26E-2 & 4.04E-5 & 2.94  & 1.63E-2 & 5.49E-6 &  3.97   \\
			\hline
		\end{tabular}		
	\end{center}  
\end{table}

\begin{figure}
	\centering
	\includegraphics[width=1.0\linewidth]{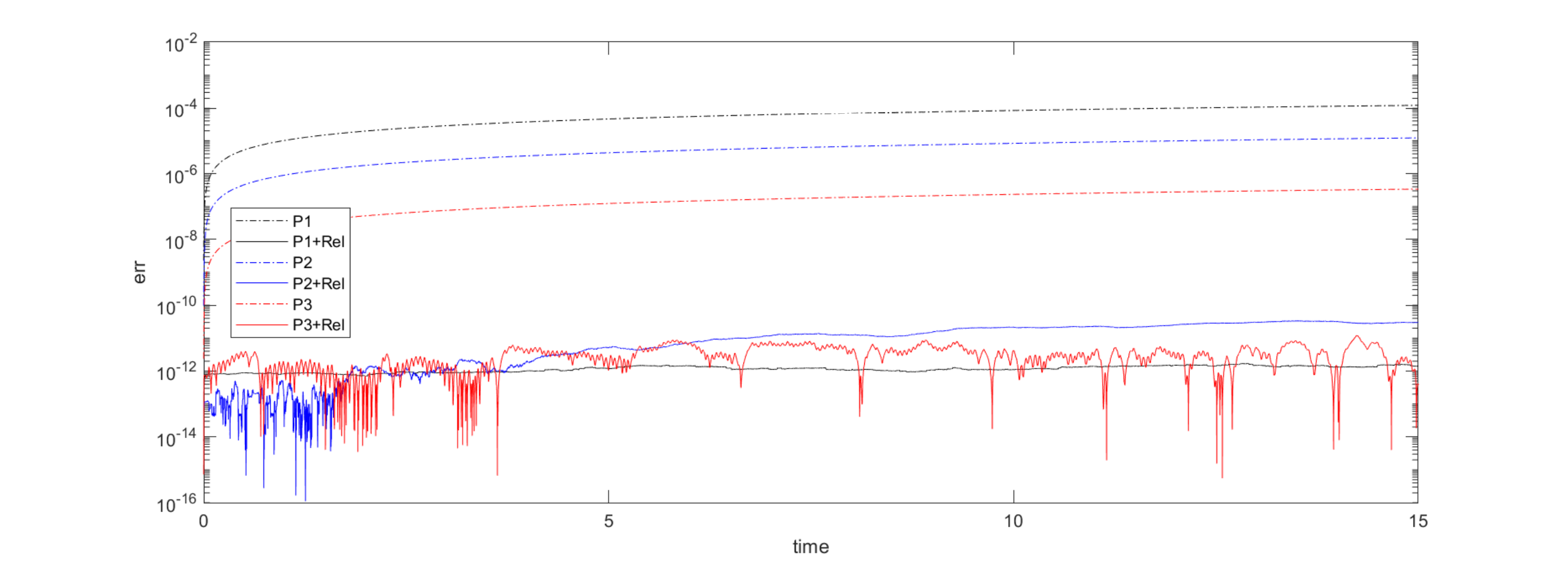}
	\caption{Total entropy error for the traveling bump test case. 
	}
	\label{fig:linear_rotation_entropy}
\end{figure}

\begin{figure}
	\centering
	\includegraphics[width=0.5\linewidth]{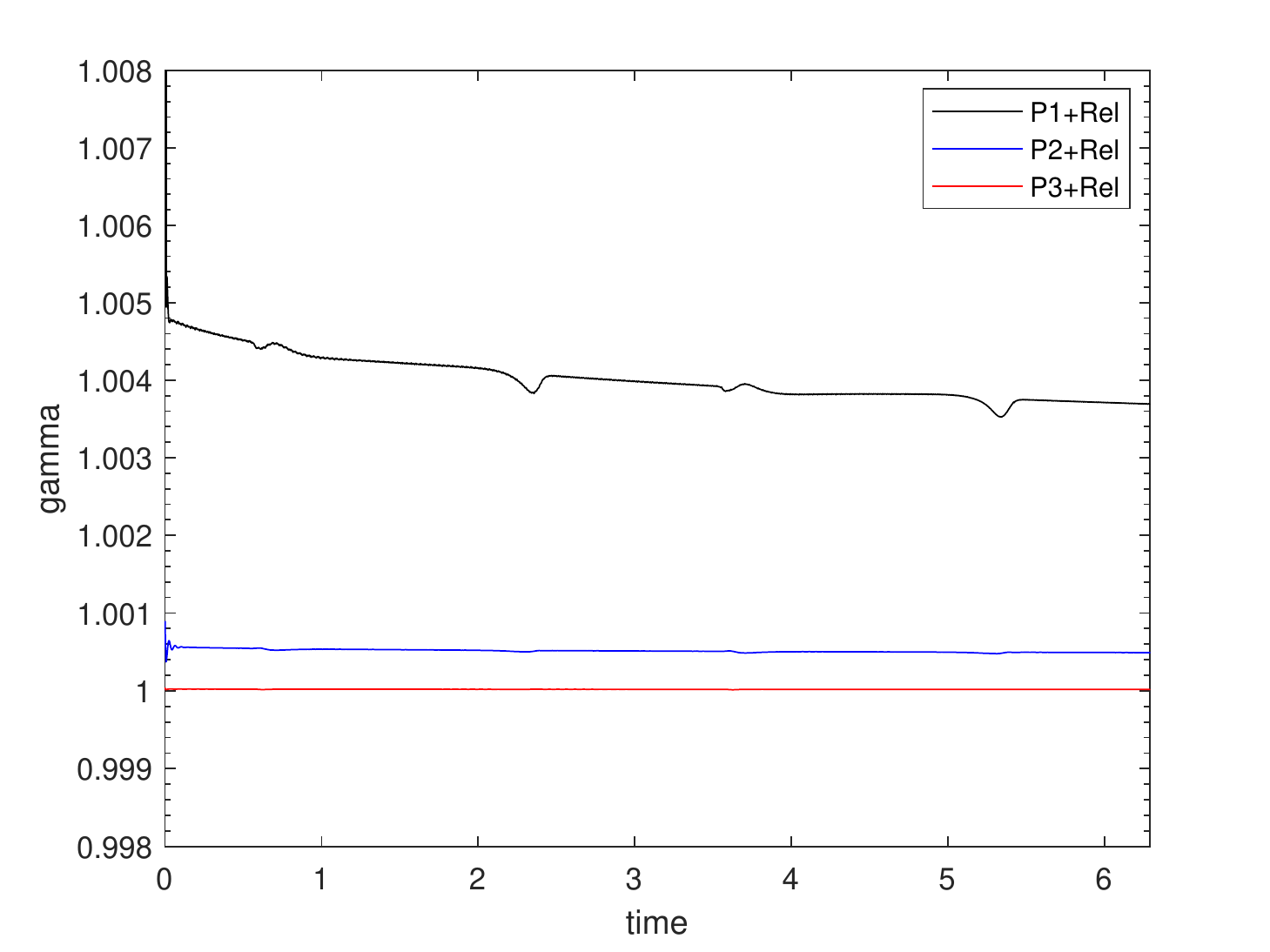}
	\caption{Values of $\gamma^n$ for the traveling bump test case. 
}
	\label{fig:linear_rotation_gamma}
\end{figure}

First of all, we check that the expected order of accuracy of the relaxed schemes is conserved. In \cref{tab.TravellingBump} we analyze the error at time $t_f=3$ with respect to the exact solution and we observe that for all elements $P_N$ we obtain order of accuracy $N+1$.
Then, we compute the error between the exact total entropy and the approximated one. In \cref{fig:linear_rotation_entropy} we compare the total entropy error for simulations with different order of polynomials, with and without relaxation, but with a similar number of degrees of freedom. In particular we run the $P_1$ scheme over a mesh of $24190$ elements ($72570$ DOFs), the $P_2$ scheme over a mesh of $11726$ elements ($70356$ DOFs), 
and the $P_3$ scheme over a mesh of $7104$ elements ($71040$ DOFs). 
We run all the simulations for time $t_f=15$.  The plot always shows the error in absolute value, but we observe in the simulations that the classical ADER-DG diffuse the total entropy. Moreover, we can see in the plot that the higher the order of accuracy the lower the error is for classical ADER-DG. Finally, we observe that the relaxation methods preserve the total entropy up to machine precision also for very long computational times.

Finally, in \cref{fig:linear_rotation_gamma} we show the values of $\gamma^n$ as a function of time for the relaxation schemes for the same meshes described above. The value is close to 1, in particular, it is much closer when the order is high.

\subsubsection{Test: rotating bump}
The second test we consider is a bump rotating around the point $(0,0)$. The initial condition is
\begin{equation}
	u_0(\vec x) = 
	\begin{cases}
		e^{1-\frac{1}{1-r^2}} & \text{if } r<1,\\
		0&\text{else},	
	\end{cases}
\end{equation}
where $r=\sqrt{x^2+(y-1.5)^2}$, the advection coefficient is $\vec a = (-y,x)^T$ and the domain is $\Omega=[-3,3]\times [-3,3]$. Homogeneous Dirichlet boundary conditions are imposed all over the boundary.

As for the previous test, we assess the order of accuracy, the error of the total entropy and the gamma values. 
The order of accuracy is computed in \cref{tab.scalar_rotation_order}  and, again, we observe the expected convergence order for all the polynomials degree.

\begin{table}
	\caption{Order of convergence for the ADER-DG schemes equipped with the relaxation technique at time $t_f=0.1$ for the $L_2$ norm error of the rotating bump.} 
	\label{tab.scalar_rotation_order}
	\begin{center} 	
		\begin{tabular}{|ccc|ccc|ccc|} 
			\hline
			\multicolumn{3}{|c|}{$P_1 +$ Rel}        &   \multicolumn{3}{|c|}{$P_2 +$ Rel}      &    \multicolumn{3}{|c|}{$P_3 +$ Rel}   \\
			$\Delta x$   &   $L_2(h-h_{ex})$   & $\mathcal{O}(L_2)$ & $\Delta x$   &   $L_2(h-h_{ex})$   & $\mathcal{O}(L_2)$ & $\Delta x$   &   $L_2(h-h_{ex})$   & $\mathcal{O}(L_2)$   \\
			\hline
			4.53E-2 & 2.93E-3 &   -  & 6.37E-2 & 8.09E-5 &  -    & 8.28E-2 & 4.28E-4 &   -    \\
			4.06E-2 & 2.35E-3 & 2.09 & 5.75E-2 & 6.16E-6 & 2.69  & 7.38E-2 & 2.91E-4 &  3.32  \\ 
			3.71E-2 & 1.97E-3 & 1.94 & 5.24E-2 & 4.77E-6 & 2.75  & 6.77E-2 & 2.12E-4 &  3.71  \\
			3.43E-3 & 1.70E-3 & 1.91 & 4.84E-2 & 3.71E-6 & 3.19  & 6.24E-2 & 1.57E-4 &  3.70  \\
			\hline
		\end{tabular}		
	\end{center}  
\end{table}

\begin{figure}
	\centering
	\includegraphics[width=1.0\linewidth]{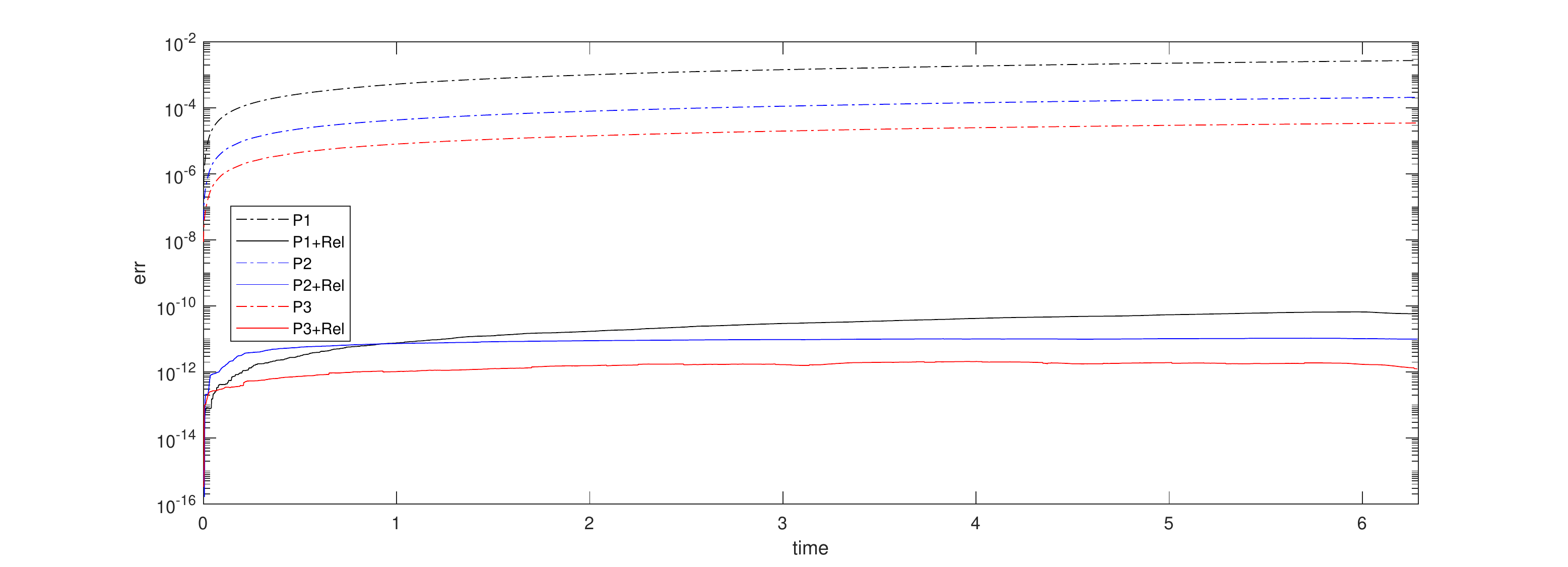}
		\caption{Total entropy error for the rotating bump test case. 
		}
	\label{fig:scalar_rotation_entropy}
\end{figure}

\begin{figure}
	\centering
	\includegraphics[width=0.5\linewidth]{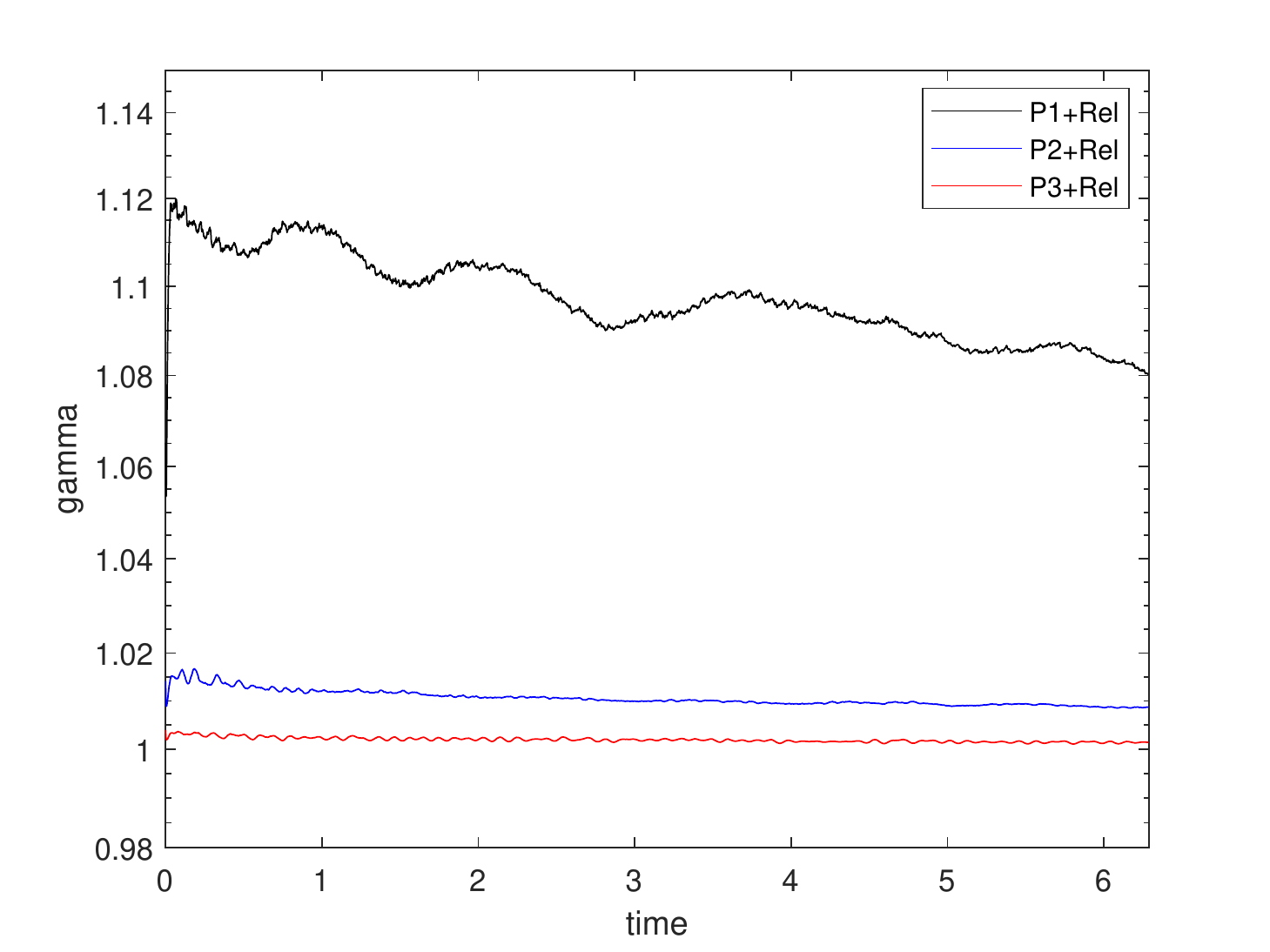}
		\caption{Values of $\gamma^n$ for the rotating bump test case. 
		}
	\label{fig:scalar_rotation_gamma}
\end{figure}

In \cref{fig:scalar_rotation_entropy} the entropy error with respect to time is depicted for the classical ADER-DG and the entropy corrected one. The three schemes considered in this plot employ a similar number of degrees of freedom:
indeed, we run the $P_1$ scheme over a mesh of $10410$ elements ($31230$ DOFs), 
the $P_2$ scheme over a mesh of $5192$ elements ($31152$ DOFs), 
and the $P_3$ scheme over a mesh of $3126$ elements ($31260$ DOFs).
In the classical ADER-DG the total entropy is preserved only up to the order of accuracy of the scheme and this error increases in time. On the other side, the relaxation ADER-DG obtains the preservation of the total entropy up to machine precision for all times.

In \cref{fig:scalar_rotation_gamma} we see the values of $\gamma^n$ as a function of time for the relaxation method. As expected, the values are close to $1$ and, as the scheme is more accurate, we obtain values closer to $1$. One must notice that for this test the values of $\gamma$ are larger with respect to the previous test. This means that the classical ADER-DG in this test has a larger entropy error with respect to the previous one, hence, a larger correction on the timestep is needed to preserve the total entropy.

\subsection{Shallow water equations}
The shallow water equations describe the motion of water for vertically averaged variables. We define with $h$ the water height and with $u$ and $v$ the vertically averaged velocity in $x$ and $y$ direction, respectively. We denote the discharge with  $\vec q:=(hu,hv)^T$ and with  $g$ the gravitational acceleration.
The shallow water equations on a flat bathymetry read
\begin{equation}
	\begin{cases}
		\partial_t h + \partial_x (hu) + \partial_y (hv) =0,\\
		\partial_t (hu) + \partial_x (hu^2 + g \frac{h^2}{2}) + \partial_y(huv)=0,\\
		\partial_t (hv) + \partial_x(huv) + \partial_y (hv^2 + g \frac{h^2}{2}) =0.
	\end{cases}
\end{equation}
The associated entropy, entropy variables and entropy fluxes are 
	\begin{align}
		\eta(\u)=h k+\frac{gh^2}{2}, \quad \v(\u) = (gh-k, u, v)^T,\quad G = \vec q(gh + k),\\
	\end{align}
	having introduced the kinetic energy $k=(u^2+v^2)/2 $. In this case we have 
	\begin{equation}
	\partial_{\u}\v =\dfrac{1}{h}\left(
	\begin{array}{ccc}
		gh+2k & -u & -v \\ 
		-u  & 1 & 0 \\ 
		-v & 0 & 1  
	\end{array}\right)\quad\text{and}\quad
	A_0=\dfrac{1}{g}\left(
	\begin{array}{ccc}
		1 & u & v \\ 
		u  & gh + u^2 & uv \\ 
		v & uv & gh + v^2  
	\end{array}\right).
	\end{equation}

\subsubsection{Moving vortex test}
For these equations, we consider some traveling vortexes whose analytical solution is available. There are many vortexes known in literature that can serve the scope. In ~\cite{ricchiuto2009stabilized} there is a $\mathcal{C}^2$ vortex and in~\cite{shuosher1} there is an isentropic Euler vortex that can be adapted to shallow water which is $\mathcal{C}^\infty$, but it is not compactly supported. In order to have a vortex that is both compactly supported and with enough smoothness, we will use one of the vortexes defined in
\cite{ricchiuto2021analytical}. In particular, we will define it on the domain $\Omega = [0,1]^2$ and such that it is $\mathcal{C}^6(\Omega)$. It is centered in $\x_c=(0.5,0.5)^T$ with radius $r_0=0.45$. Its amplitude is $\Delta h =0.1$, its frequency is $\omega=\pi /r_0$ and its intensity is $\Gamma = \frac{12 \pi \sqrt{g \Delta h } }{r_0 \sqrt{315 \pi^2-2048}}$.
Then, we define the reference water level $h_c=1$ and the background velocity $u_c=1$, $v_c=0$. Next, we introduce the co-moving coordinates $\mathcal{I}(\mathbf{x},t) = \mathbf{x} - \mathbf{x_c} - (u_c t,v_c t)^T$ and we define the distance from the vertex center as $\mathcal{R}( \mathbf{x},t) = \|   \mathcal{I}(\mathbf{x},t) \|$.

\begin{figure}
	\centering
	\includegraphics[trim= 5 5 5 5,clip,width=0.33\linewidth]{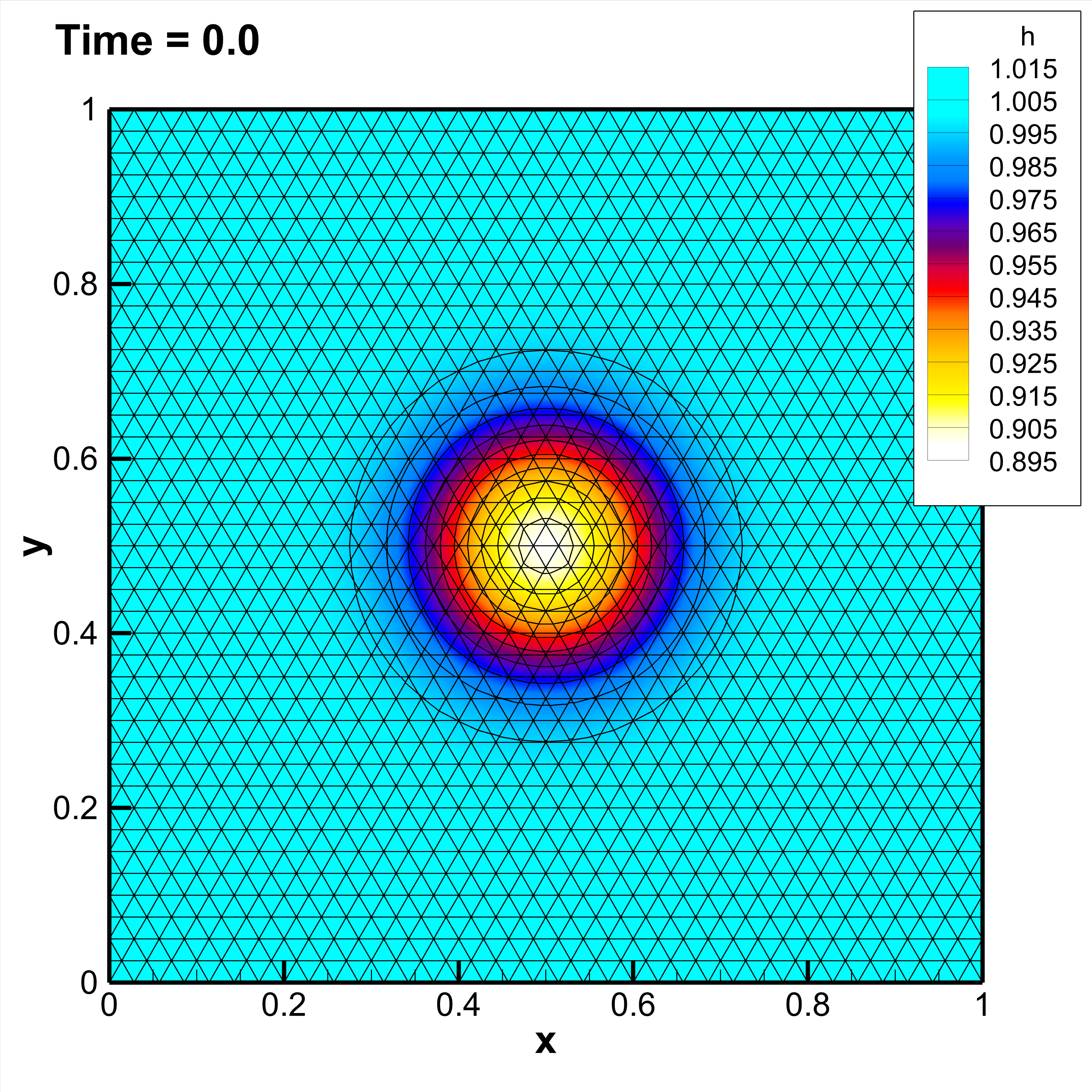}
	\includegraphics[trim= 5 5 5 5,clip,width=0.33\linewidth]{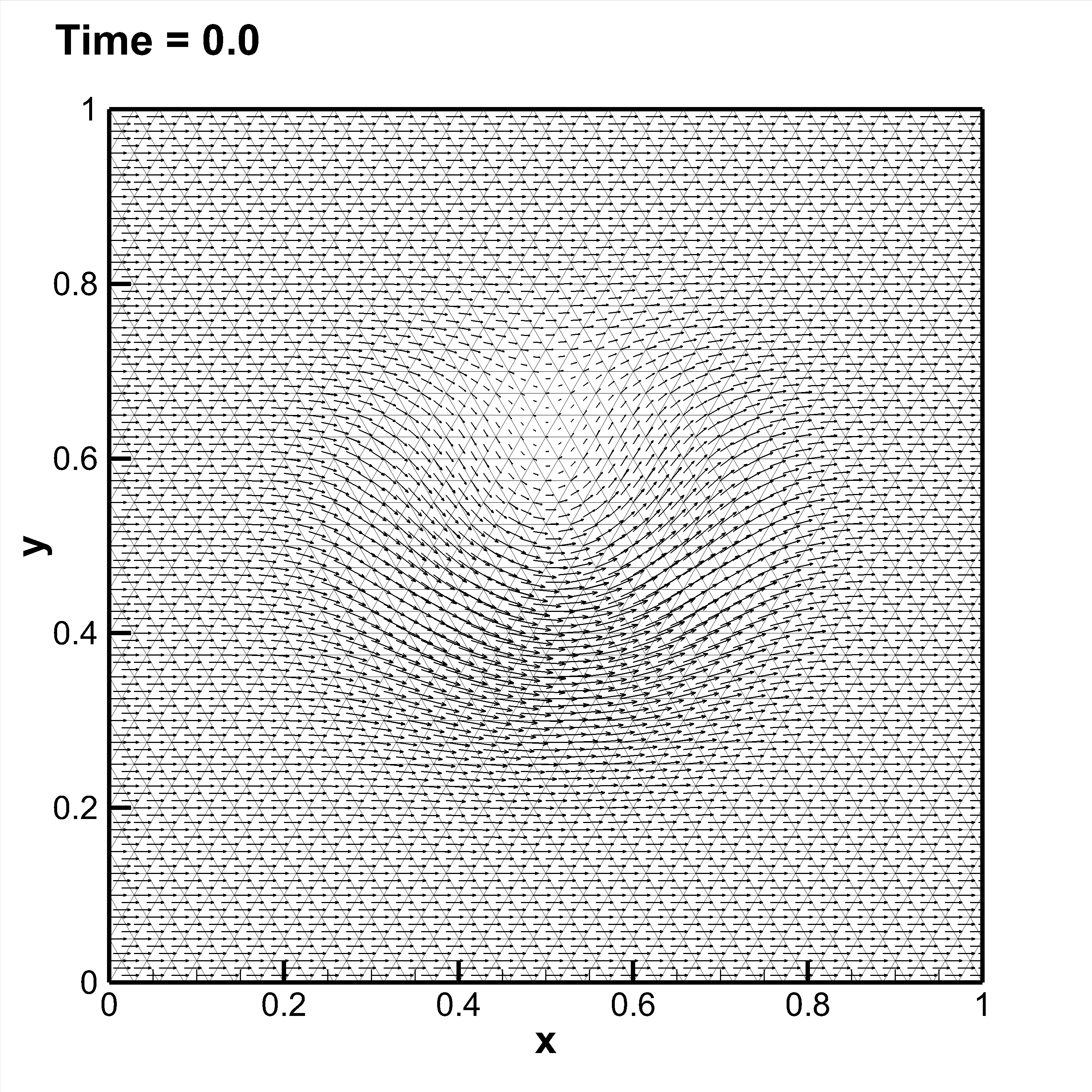}
	\includegraphics[trim= 5 5 5 5,clip,width=0.33\linewidth]{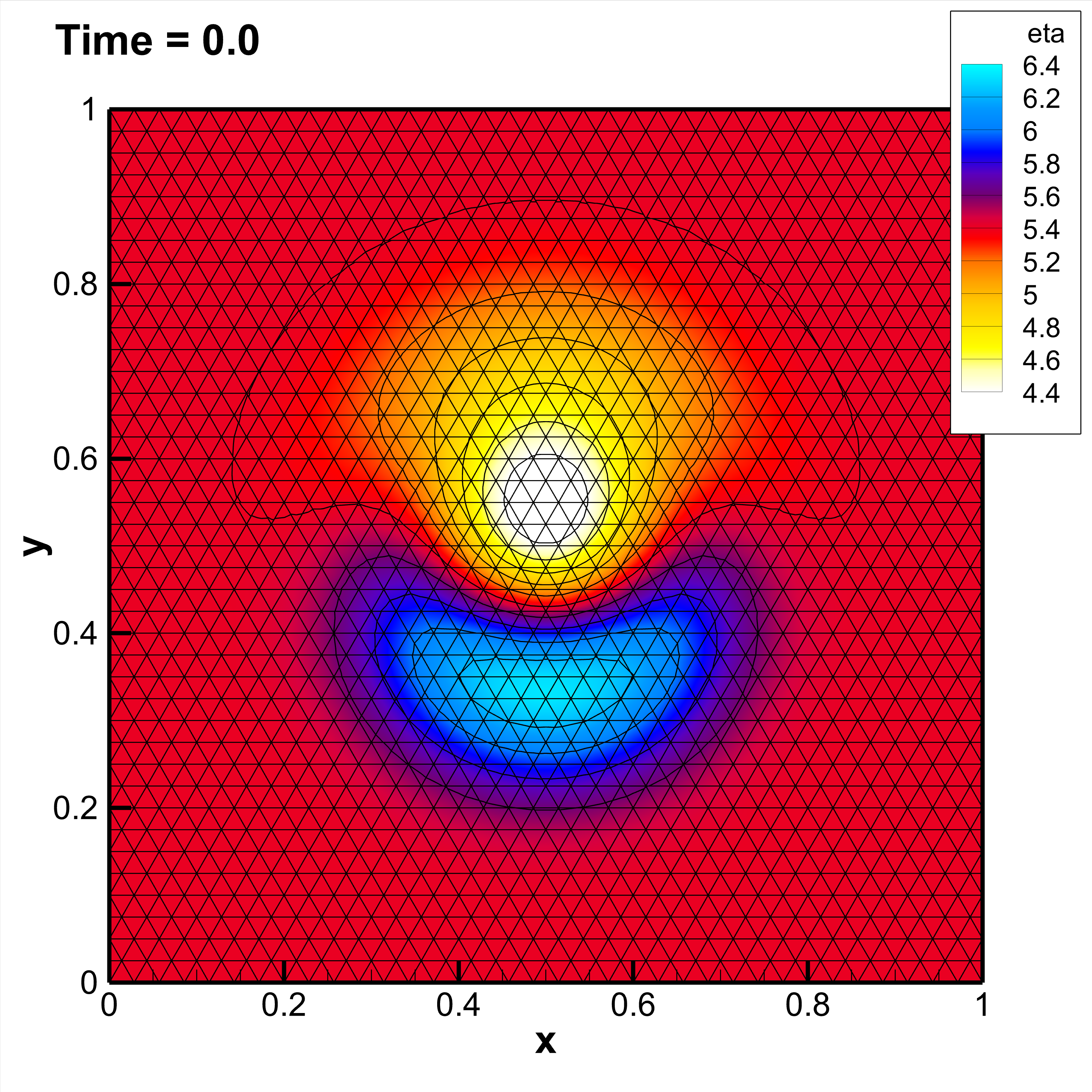} \\
	\includegraphics[trim= 5 5 5 5,clip,width=0.33\linewidth]{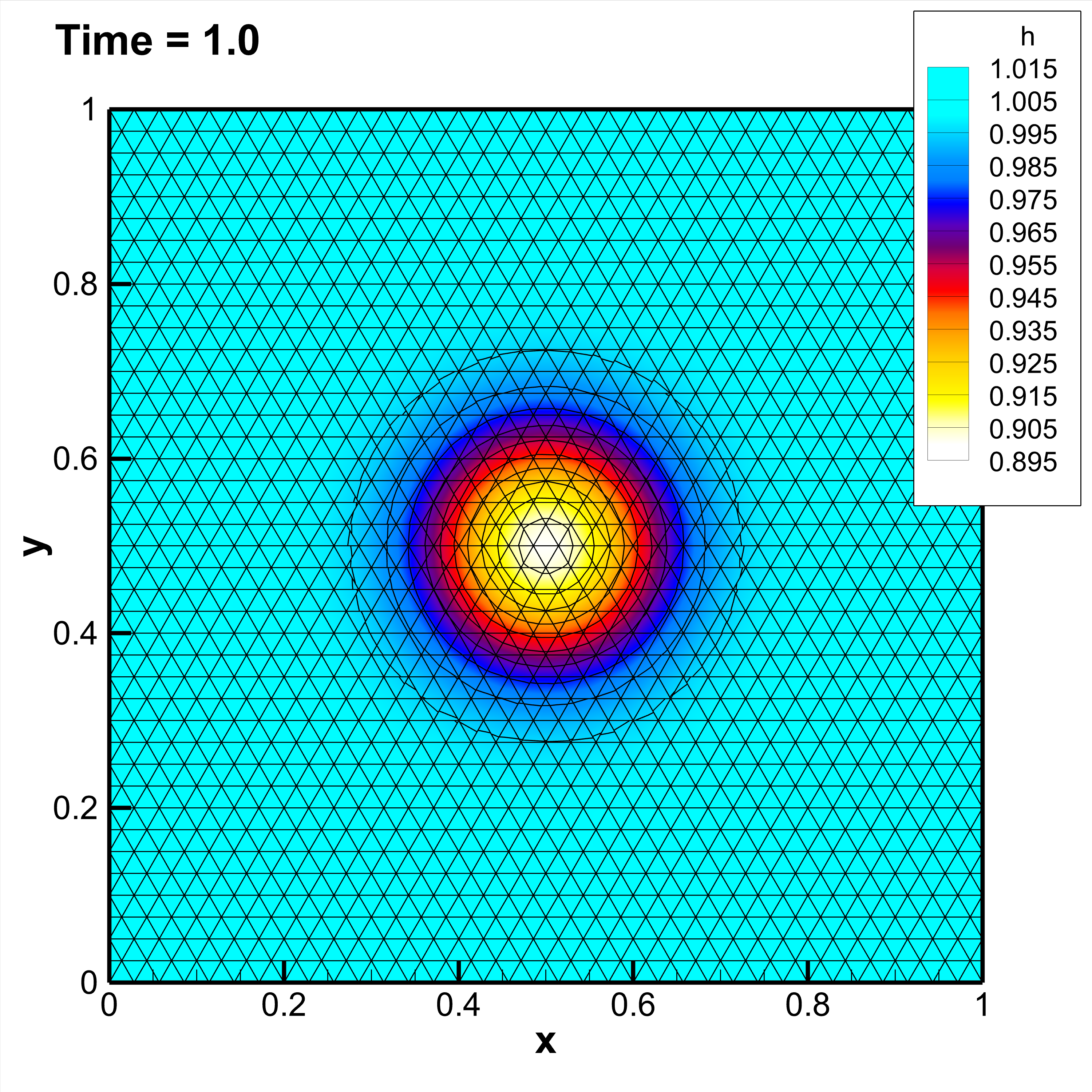}
	\includegraphics[trim= 5 5 5 5,clip,width=0.33\linewidth]{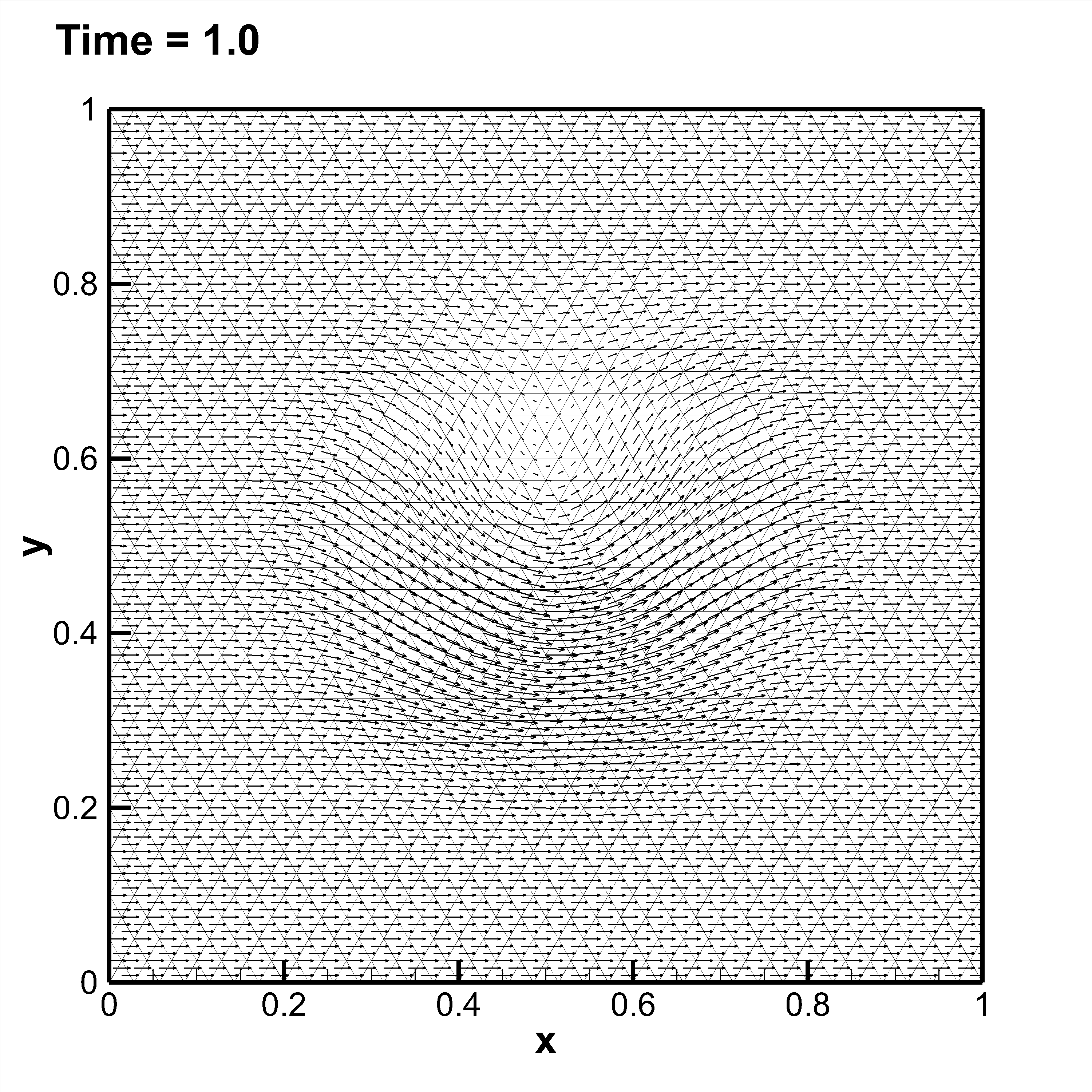}
	\includegraphics[trim= 5 5 5 5,clip,width=0.33\linewidth]{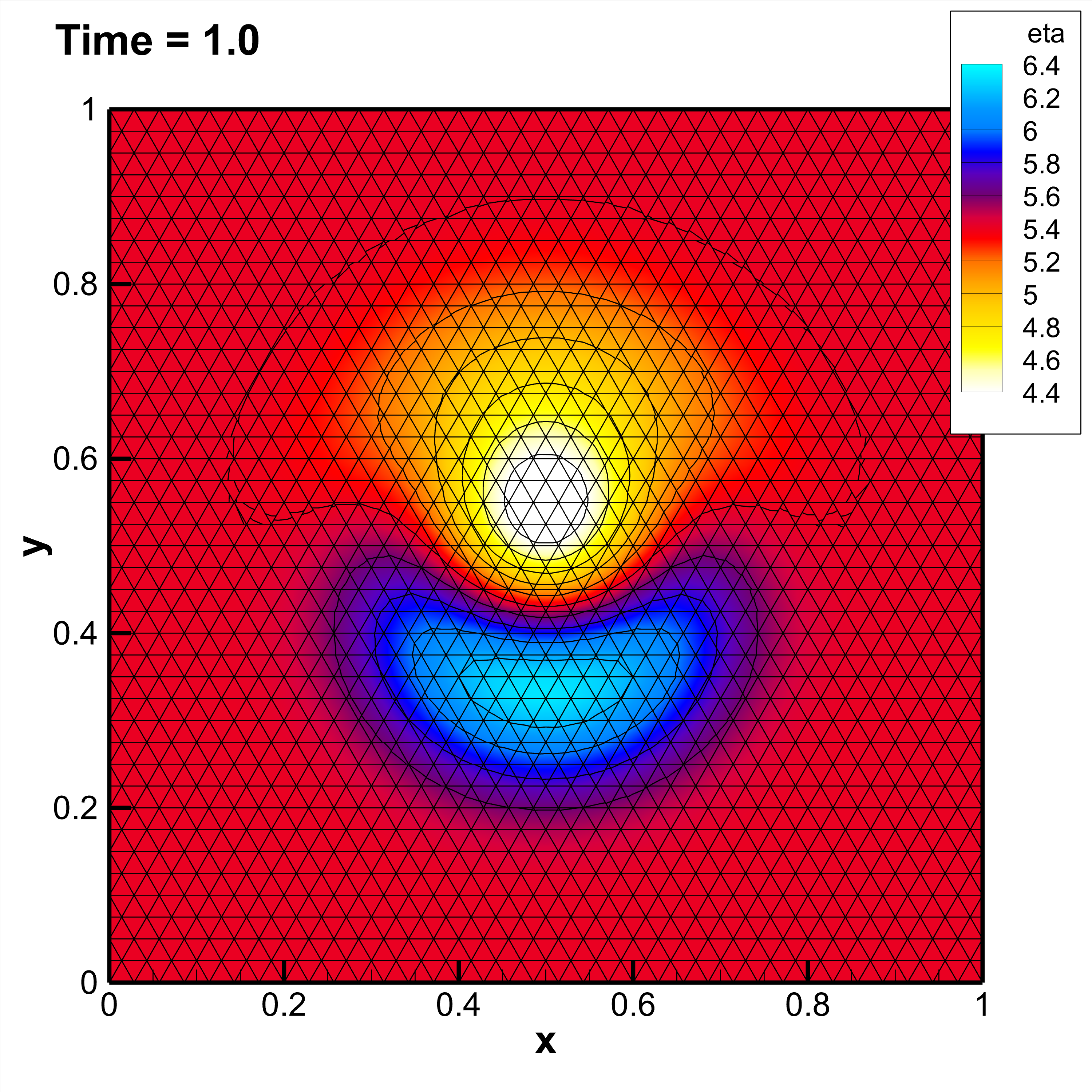}
	\caption{SW vortex test case. 
		In the figure we report the profile of the density $\rho$ (left), the velocity field (middle) and the entropy $\eta$ (right) at time $t=0$ (top row) and $t=1$ (bottom row) obtained with a DG scheme of order $2$ with relaxation on a coarse mesh made of $2840$ elements.
	}
	\label{fig:SWvortexplot}
\end{figure}

The analytical solution is defined as
\begin{align}\label{eq_sol2_SW2}
	\begin{pmatrix}
		h(x,t)\\u(x,t)\\v(x,t)
	\end{pmatrix}=
	\begin{cases}
		\begin{pmatrix}
			h_c + \frac{1}{g} \frac{\Gamma^2}{\omega^2} \cdot \left( \lambda(\omega  \mathcal{R}( \mathbf{x},t) ) - \lambda (\pi) \right) ,  \\
			u_c + \Gamma(1+\cos (\omega \mathcal{R}( \mathbf{x},t)))^2 \cdot (- \mathcal{I}(\mathbf{x},t)_y),  \\
			v_c + \Gamma(1+\cos (\omega \mathcal{R}( \mathbf{x},t)))^2 \cdot ( \mathcal{I}(\mathbf{x},t)_x),
		\end{pmatrix}, &\mbox{if } \omega  \mathcal{R}( \mathbf{x},t) \leq \pi,\\
		\begin{pmatrix}
			h_c & u_c & v_c
		\end{pmatrix}^T, &\mbox{else,}
	\end{cases}
\end{align}    
with 
\begin{align*}
	\lambda(r) = &\frac{20\cos(r)}{3} + \frac{27\cos(r)^2}{16} + \frac{4\cos(r)^3}{9} + \frac{\cos(r)^4}{16} + \frac{20r\sin(r)}{3} \\
	&+ \frac{35r^2}{16} + \frac{27r\cos(r)\sin(r)}{8} + \frac{4r\cos(r)^2 \sin(r)}{3} + \frac{r\cos(r)^3 \sin(r)}{4}.
\end{align*}
We consider periodic boundary conditions and final time $t_f=1$. A plot of the initial solution and of the solution at $t_f=1$ simulated with the relaxation ADER-DG $P_2$ is given in Figure~\ref{fig:SWvortexplot}.

For this test we study the order of accuracy of the schemes for polynomials of degrees 1, 2 and 3. In \cref{tab.SWvortex} we list the errors $L_2(h-h_{ex})$ and the experimental order of accuracy $\mathcal O(L_2)$. The order of accuracy found is the expected one, i.e., 2, 3, and 4. It is clear that the error and the accuracy of the scheme is not modified by the entropy correction and relaxation techniques.

\begin{table}
	\caption{Order of convergence for the ADER-DG schemes equipped with the relaxation technique at time $t_f=1$ for the $L_2$ norm error of the water height $h$ in the shallow water moving vortex test case.} 
	\label{tab.SWvortex}
	\begin{center} 	
		\begin{tabular}{|ccc|ccc|ccc|} 
			\hline
			\multicolumn{3}{|c|}{$P_1 +$ Rel}        &   \multicolumn{3}{|c|}{$P_2 +$ Rel}      &    \multicolumn{3}{|c|}{$P_3 +$ Rel}   \\
			$\Delta x$   &   $L_2(h-h_{ex})$   & $\mathcal{O}(L_2)$ & $\Delta x$   &   $L_2(h-h_{ex})$   & $\mathcal{O}(L_2)$ & $\Delta x$   &   $L_2(h-h_{ex})$   & $\mathcal{O}(L_2)$   \\
			\hline
			1.64E-2 & 9.37E-5 &   -  & 2.36E-2 & 1.02E-5 &  -    & 2.98E-2 & 8.41E-7 &   -    \\
			1.27E-2 & 5.53E-5 & 2.07 & 1.83E-2 & 4.93E-6 & 2.89  & 2.36E-2 & 3.23E-7 &  4.10  \\ 
			1.01E-2 & 3.42E-5 & 2.04 & 1.43E-2 & 2.40E-6 & 2.91  & 1.84E-2 & 1.17E-7 &  4.07  \\
			7.88E-3 & 2.09E-5 & 2.02 & 1.10E-2 & 1.11E-6 & 2.93  & 1.43E-2 & 4.28E-8 &  4.05  \\
			\hline
		\end{tabular}		
	\end{center}  
\end{table}

\begin{figure}
	\centering
\includegraphics[width=1.0\linewidth]{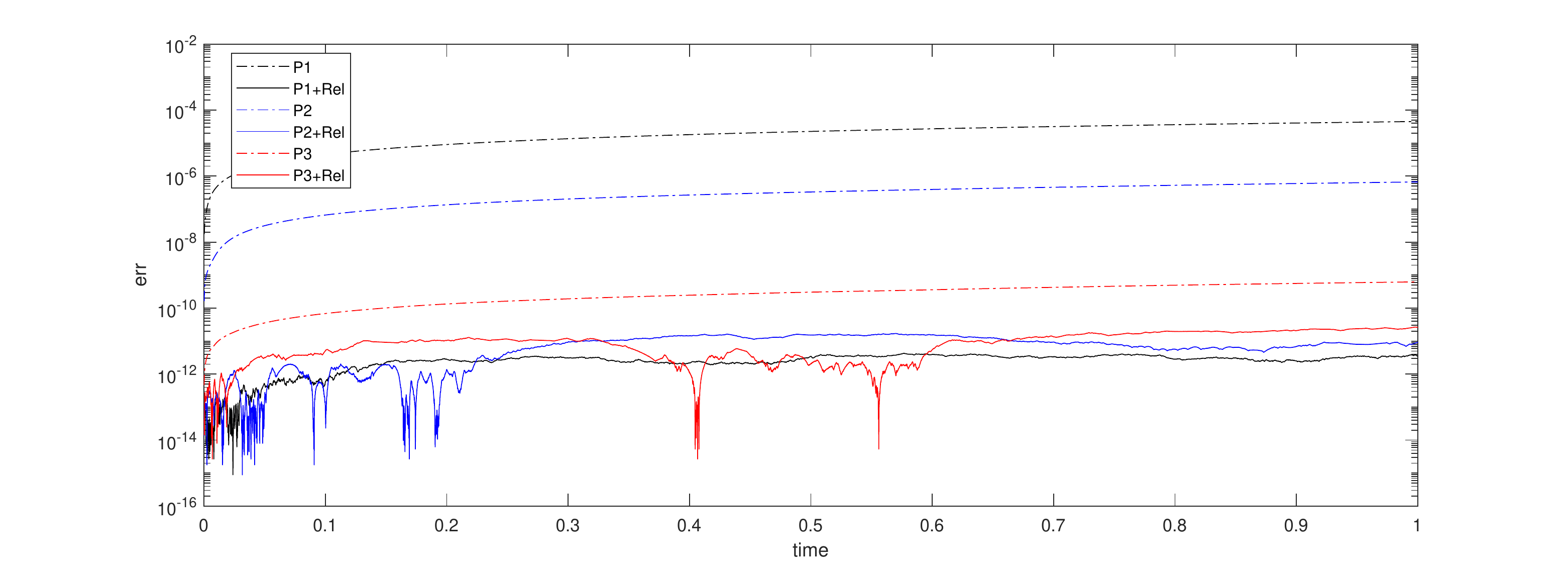}
	\caption{Global entropy error for the shallow water moving vortex test case. 
	}
	\label{fig:SWvortexentropy}
\end{figure}

To compare the results with the classical ADER-DG method, we consider different meshes: for $\mathbb P^1$ $7590$ elements ($22770$ DOFs), for $\mathbb P^2$ $3726$ elements ($22356$ DOFs) and  for $\mathbb P^3$ $2268$ elements ($22680$ DOFs), in order to have a comparable number of DOFs.
In \cref{fig:SWvortexentropy} we plot the total entropy error as a function of time for these schemes. For the classical ADER-DG methods the global entropy is dissipated in time and the lower the order the larger is the entropy error. Conversely, for the relaxation version of the schemes, the total entropy is conserved up to machine precision independently on the order of accuracy of the scheme.

\begin{figure}
	\centering
	\includegraphics[width=0.5\linewidth]{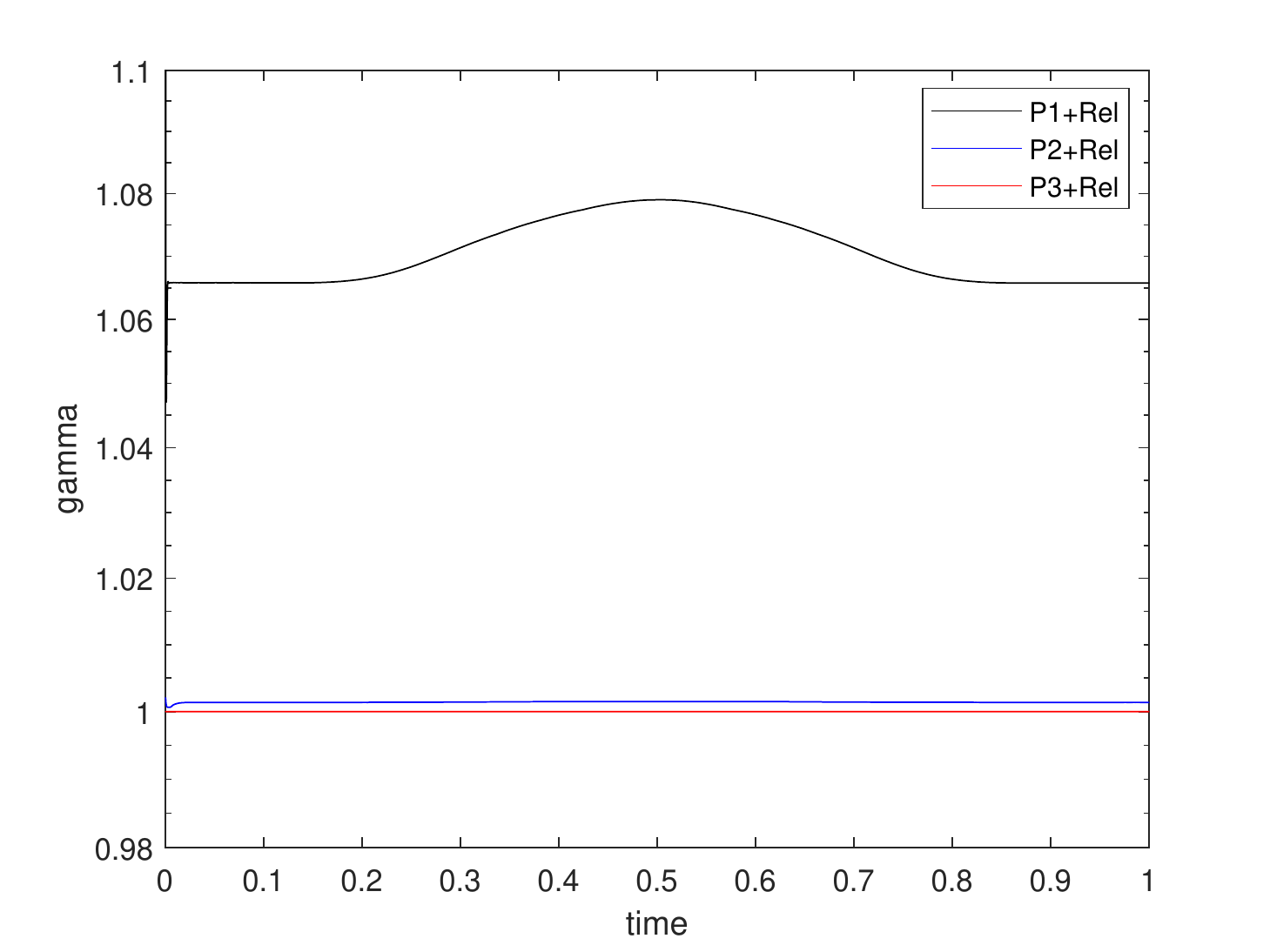}
	\caption{Values of $\gamma^n$ for the shallow water vortex test case. 
	}
	\label{fig:SWvortexgamma}
\end{figure}

In \cref{fig:SWvortexgamma} we can observe the values of $\gamma^n$ in time for the three relaxation versions of the methods. We observe that for $\mathbb P^1$ $\gamma^n=1+\mathcal{O}(10^{-1})$, while for higher orders the values of $\gamma^n$ are much closer to 1, as the original scheme is already close enough to the entropy conservative solution.

\subsection{Euler equations of gas dynamics}

Euler's equations describe the dynamic of a gas in $\rho$ the density of the gas, $\vec u=(u,\,v)$ the velocity vector and $E$ the total energy. Through these variables, we can define other quantities to shorten the notation of the equations: $k=(u^2+v^2)/2$ the specific kinetic energy, $e:=E-\rho k$ the internal energy and $p$ the pressure linked to the conservative variables by the equation of state (EOS)
\begin{equation}
	E=\frac{p}{\gamma -1} +\rho k,
\end{equation}
where $\gamma$ is the adiabatic constant.
Then, Euler's equations for a perfect gas read 
\begin{equation}
\begin{cases}
\partial_t \rho + \partial_x (\rho u) +\partial_y(\rho v ) =0, \\
\partial_t (\rho u) + \partial_x (\rho u^2 +p) + \partial_y(\rho u v) =0,\\
\partial_t (\rho v) + \partial_x (\rho uv) + \partial_y(\rho v^2 + p) =0,\\
\partial_t (E) + \partial_x (u (E+p) )+ \partial_y (v (E+p) )=0.
\end{cases}
\end{equation}
Denoting with $s=p\rho^{-\gamma}$, we consider the entropy 
\begin{equation}
\eta(\u) = -\frac{\gamma +1}{\gamma -1}\rho s^{\frac{1}{\gamma +1}}.
\end{equation}
Let us denote also the momentum as $\vec{m}^T := (\rho u, \rho v)$, then, the entropy variables  and entropy fluxes are
\begin{equation}
\v(\u)= - (\rho p)^{-\frac{\gamma}{\gamma+1}} \begin{pmatrix}
E\\ -\vec{m} \\ \rho
\end{pmatrix},\quad
G = -\vec m \frac{\gamma+1}{\gamma-1}(p\rho^{-\gamma})^{\frac{1}{\gamma+1}}.
\end{equation}
The Hessian of the flux is
\begin{equation}
\partial_{\u}\v = \frac{\gamma (\gamma-1)}{\gamma +1}(\rho p)^{-\frac{2\gamma+1}{\gamma+1}} \begin{pmatrix}
E& -\vec{m} & \rho
\end{pmatrix}  \otimes \begin{pmatrix}
E\\ -\vec{m} \\ \rho
\end{pmatrix}-(\rho p)^{-\frac{\gamma}{\gamma+1}}\begin{pmatrix}
0&0&0&1\\
0&-1&0&0\\
0&0&-1&0\\
1&0&0&0
\end{pmatrix},
\end{equation}
and $A_0$, the inverse of $\partial_{\u}\v$, is
\begin{equation}
A_0=\gamma \rho (\rho p)^{-\frac{1}{\gamma+1}} \begin{pmatrix}
\rho& \rho u & \rho  v& \rho k +\frac{p}{\gamma(\gamma-1)}\\
\rho u & E- \frac{p}{\gamma(\gamma-1)}+\frac{\rho (u^2-v^2)}{2} &\rho uv & uE\\
\rho v &\rho uv & E- \frac{p}{\gamma(\gamma-1)}-\frac{\rho (u^2-v^2)}{2}  & vE\\
\rho k + \frac{p}{\gamma(\gamma-1)}& Eu &E v & \frac{E^2}{\rho}		
\end{pmatrix}.
\end{equation}

\subsubsection{Shu vortex}
\begin{figure}
	\centering
	\includegraphics[width=0.33\linewidth]{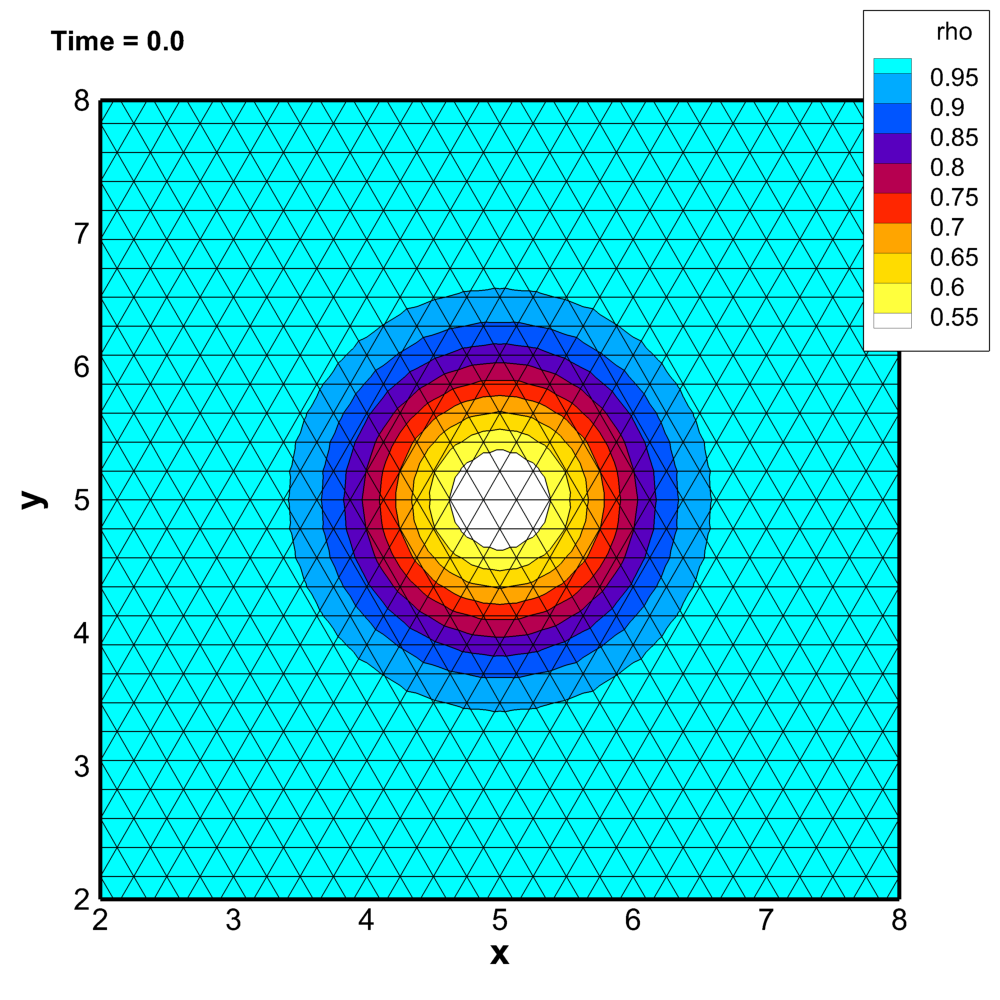}
	\includegraphics[width=0.33\linewidth]{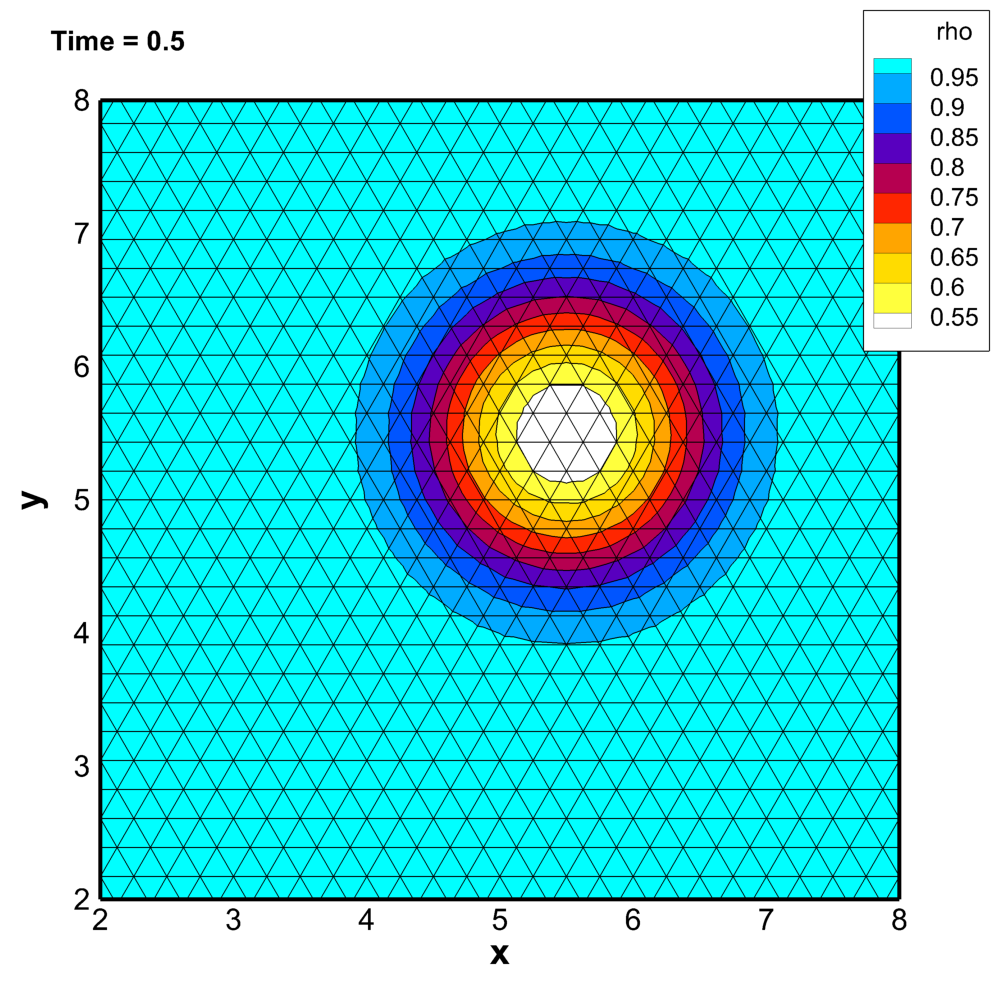}
	\includegraphics[width=0.33\linewidth]{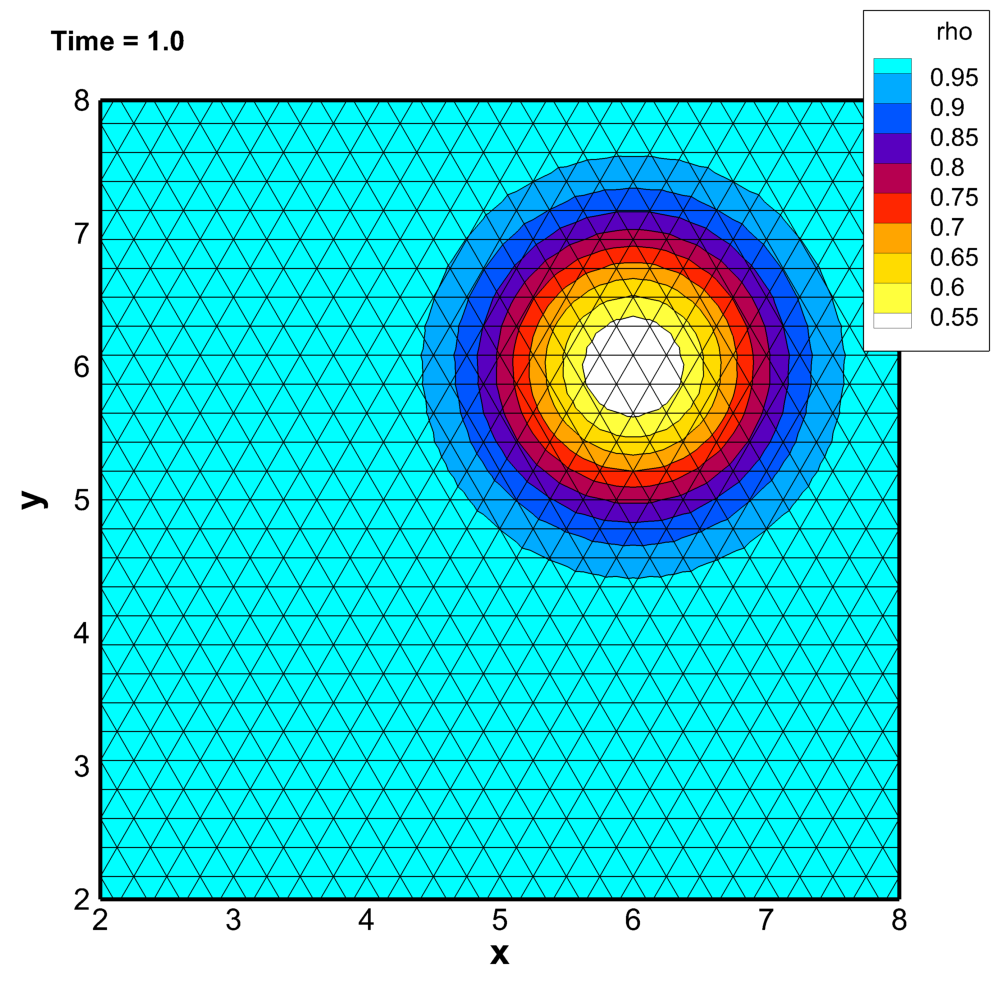}
	\includegraphics[width=0.33\linewidth]{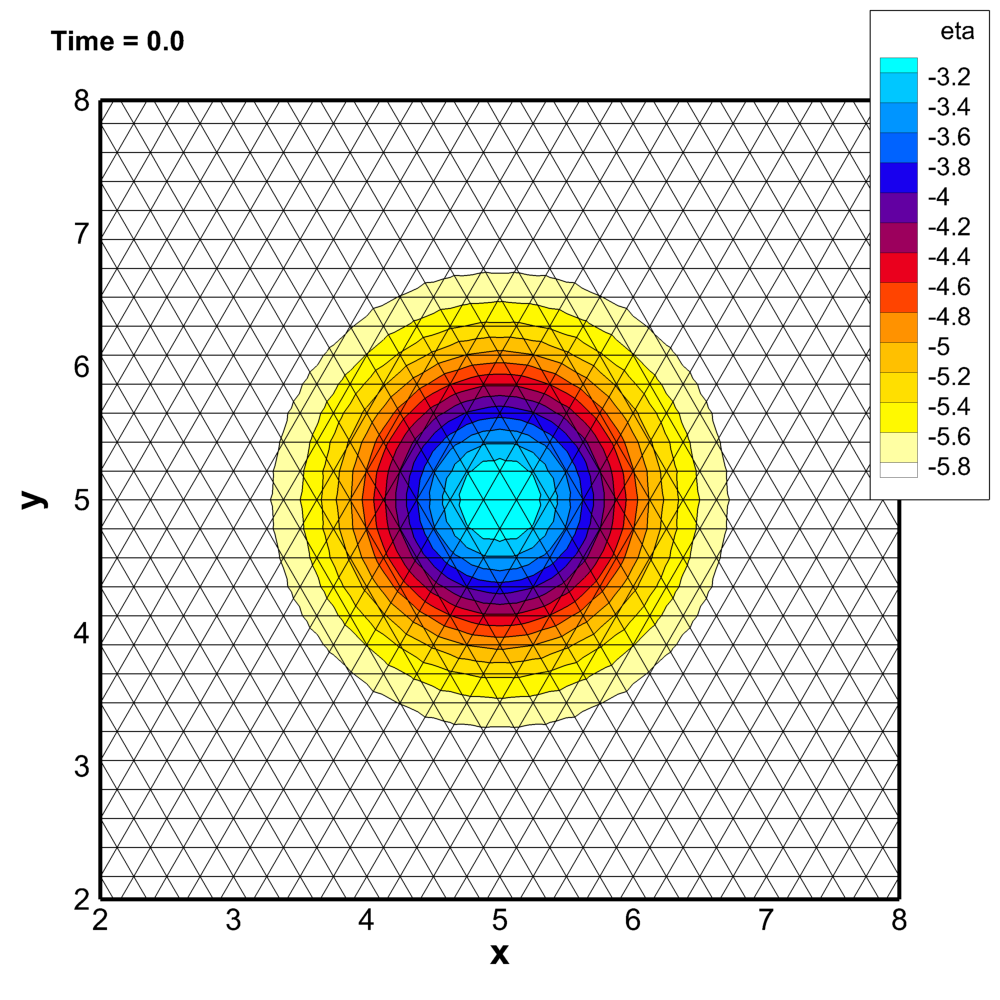}
	\includegraphics[width=0.33\linewidth]{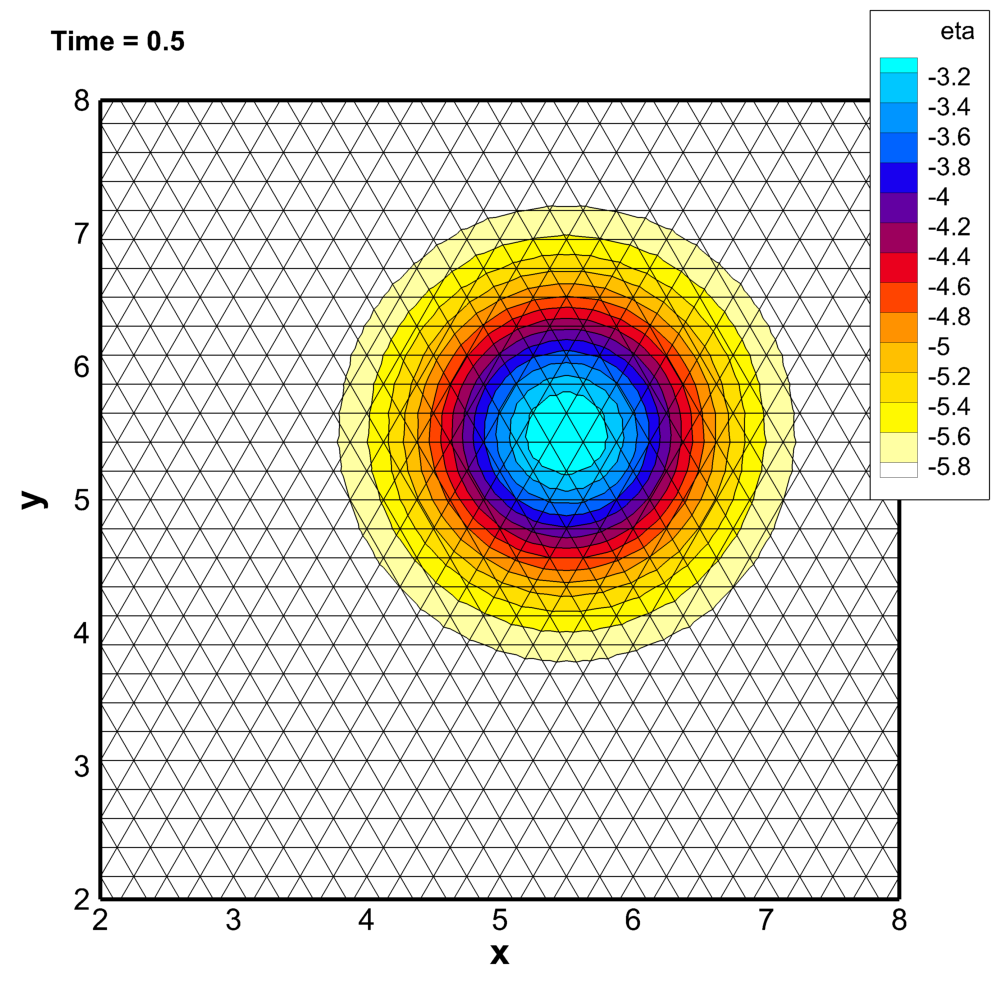}
	\includegraphics[width=0.33\linewidth]{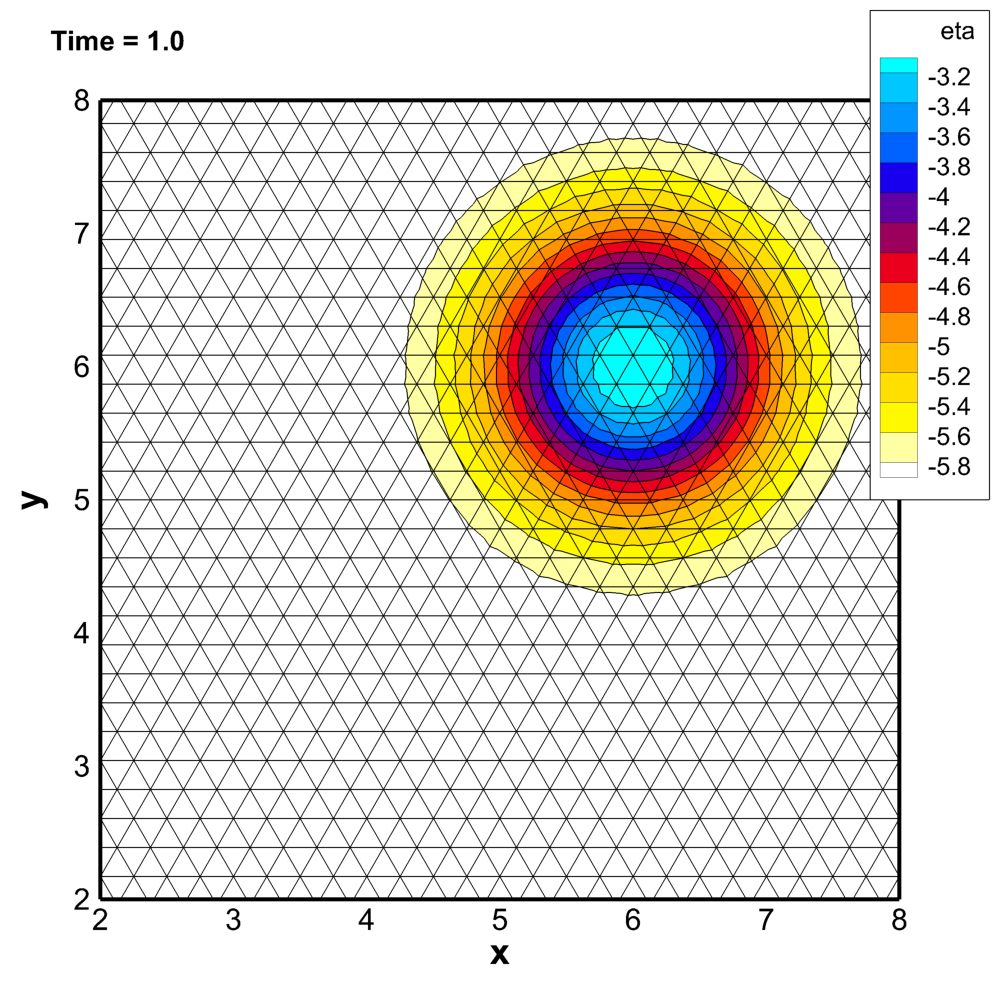}
	\caption{Shu vortex test case. 
		In the figure we report the profile of the density $\rho$ (top row) and the entropy $\eta$ (bottom row) at time $t=0, 0.5, 1$ obtained with a DG scheme of order $4$ with relaxation on a very coarse mesh made of $3726$ elements.
	}
	\label{fig:eulershuplot}
\end{figure}

We consider again a moving vortex to assess the accuracy of the scheme and to compare the entropy production of the proposed schemes. We consider a Shu vortex~\cite{shuosher1} at the beginning centered in $(5,5)$ with periodic boundary conditions.  
We plot some simulations at different times of the variable $\rho$ and of $\eta(\u)$ in \cref{fig:eulershuplot} for the DG of order 4 on a very coarse mesh.
We test first of all the order of accuracy of the proposed schemes. In \cref{tab.EulerShu} we observe that all the methods behave with the expected accuracy. 

\begin{table}
		\caption{Order of convergence for the ADER-DG schemes equipped with the relaxation technique at time $t_f=1$ for the $L_2$ norm of the density variable $\rho$ in the Shu vortex test case. } 
		\label{tab.EulerShu}
		\begin{center} 	
			\begin{tabular}{|ccc|ccc|ccc|} 
				\hline
				\multicolumn{3}{|c|}{$P_1 +$ Rel}        &   \multicolumn{3}{|c|}{$P_2 +$ Rel}      &    \multicolumn{3}{|c|}{$P_3 +$ Rel}   \\
				$\Delta x$   &   $L_2(\rho-\rho_{ex})$   & $\mathcal{O}(L_2)$ & $\Delta x$   &   $L_2(\rho-\rho_{ex})$   & $\mathcal{O}(L_2)$ & $\Delta x$   &   $L_2(\rho-\rho_{ex})$   & $\mathcal{O}(L_2)$   \\
				\hline
				1.01E-1 & 1.48E-3 &   -  & 1.43E-1 & 3.10E-4 &  -    & 1.84E-1 & 1.81E-5 &   -    \\
				7.88E-2 & 9.02E-4 & 2.03 & 1.10E-1 & 1.54E-4 & 2.68  & 1.43E-1 & 6.18E-6 &  4.35  \\ 
				5.63E-2 & 4.58E-4 & 2.02 & 8.09E-2 & 6.65E-5 & 2.71  & 1.04E-1 & 1.60E-6 &  4.20  \\
				4.11E-2 & 2.43E-4 & 2.01 & 5.82E-2 & 2.69E-5 & 2.74  & 7.55E-2 & 4.37E-7 &  4.06   \\
				\hline
			\end{tabular}		
		\end{center}  
\end{table}

\begin{figure}
	\centering
	\includegraphics[width=1.0\linewidth]{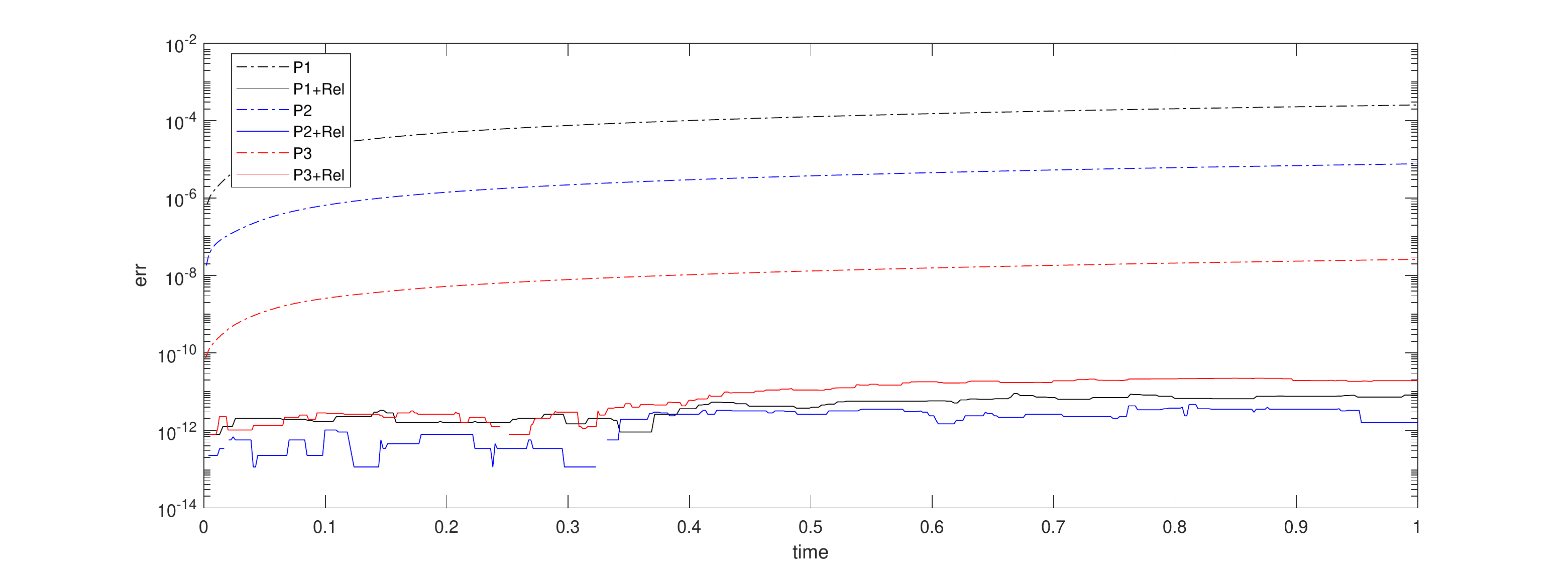}
	\caption{Entropy difference for the Shu vortex test case. 
	}
	\label{fig:eulershuentropy}
\end{figure}

Even if the error of the classical ADER-DG and the relaxation version are very close, there is a huge difference between the schemes on the total entropy. In \cref{fig:eulershuentropy} we plot the total entropy error for different simulations. We considered both classical ADER-DG and relaxation schemes for different orders.
The meshes are for $\mathbb P_1$  of $12348$ elements ($37044$ DOFs), for $\mathbb P_2$ scheme  of $6300$ elements ($37800$ DOFs), 
and for $\mathbb P_3$  of $3726$ elements ($37260$ DOFs), so that the DOFs are comparable.
As in the previous tests, we have that classical ADER-DG dissipates the total entropy, proportionally to the accuracy of the scheme, while the relaxation version preserves it up to machine precision.

\begin{figure}
	\centering
	\includegraphics[width=0.5\linewidth]{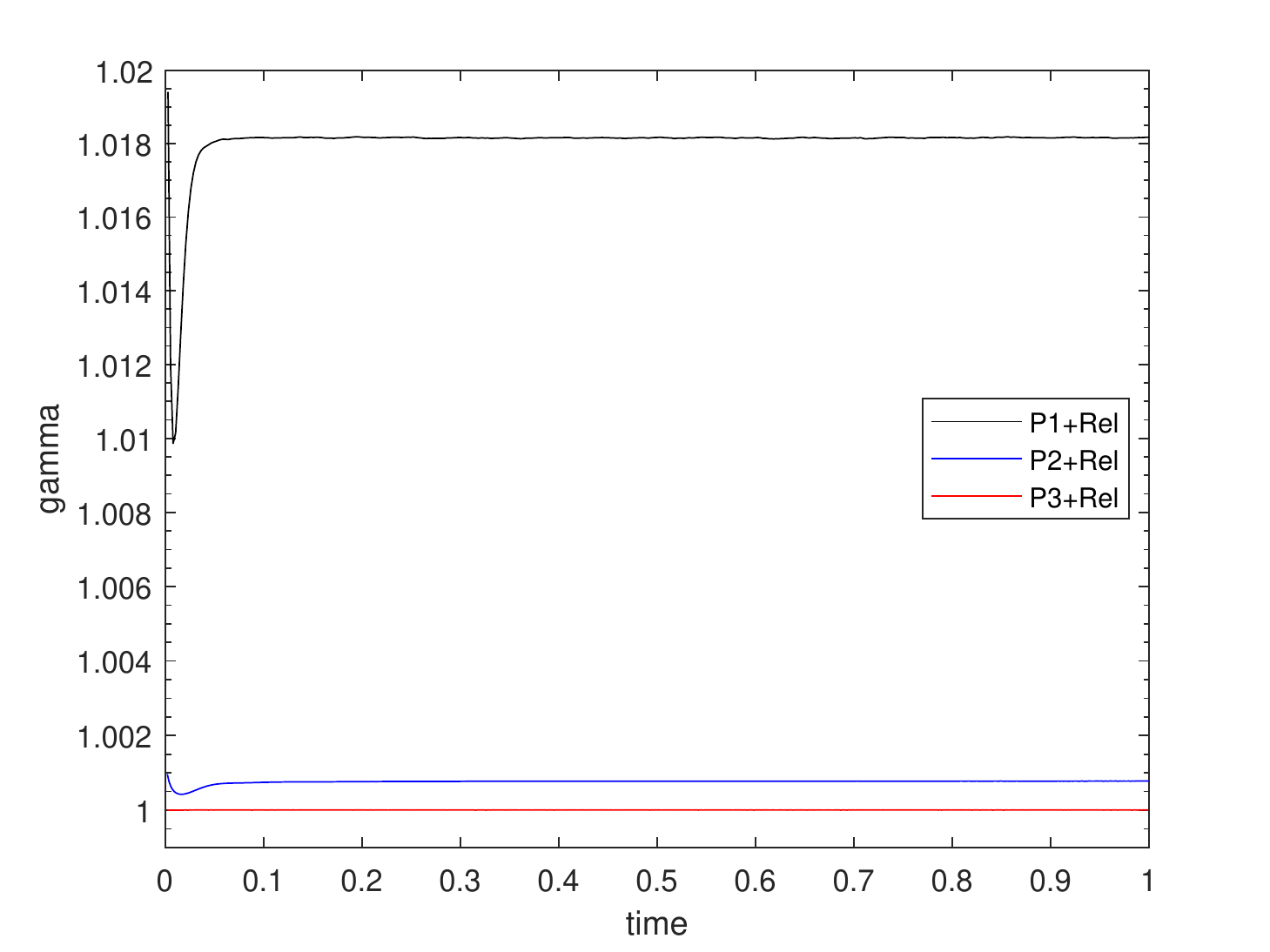}
	\caption{Values of $\gamma^n$ for the Shu vortex test case. 
}
	\label{fig:eulershugamma}
\end{figure}

\subsubsection{Moving contact discontinuity}
Now, we perform a simulation involving a moving contact discontinuity. We define a one dimensional Riemann problem on the rectangle $\Omega = [-1,1] \times [0,1]$ and the initial conditions
\begin{equation}
	\u =\begin{cases}
		\u_L &  x<0,\\
		\u_R &  x\geq 0,
	\end{cases}= \begin{cases}
		(\rho_L,\vec u_L, p_L) &  x<0,\\
		(\rho_R,\vec u_R, p_R) &  x\geq 0.
	\end{cases}
\end{equation}
We use the following values of the left and right state
\begin{equation}
		(\rho_L, u_L, v_L, p_L) =(1.5,1,0,1), \qquad (\rho_R,u_R, v_R, p_R) =(1,1,0,1).
\end{equation}
We have also to remark that we slightly smoothen the initial discontinuity (to avoid the need of the limiter) according to~\cite{tavelli2014high}
\begin{equation}
	\u = \frac{1}{2}\left( \u_R+\u_L\right) + \frac{1}{2} \left (\u_R - \u_L \right)\text{erf}\left( \frac{x}{2 \bar{h}}\right),
\end{equation}
where $\bar{h}$ is the average mesh size.
The solution of this Riemann problem is a contact discontinuity, hence, the entropy is preserved across it in time.
We consider periodic boundary conditions on the top and bottom boundaries and inflow/outflow at the left and right boundaries. We can compute analytically the amount of total entropy that leave the domain at the right boundary to compute the error with respect to the analytical solution. 

\begin{figure}
	\centering
	\includegraphics[width=0.49\linewidth]{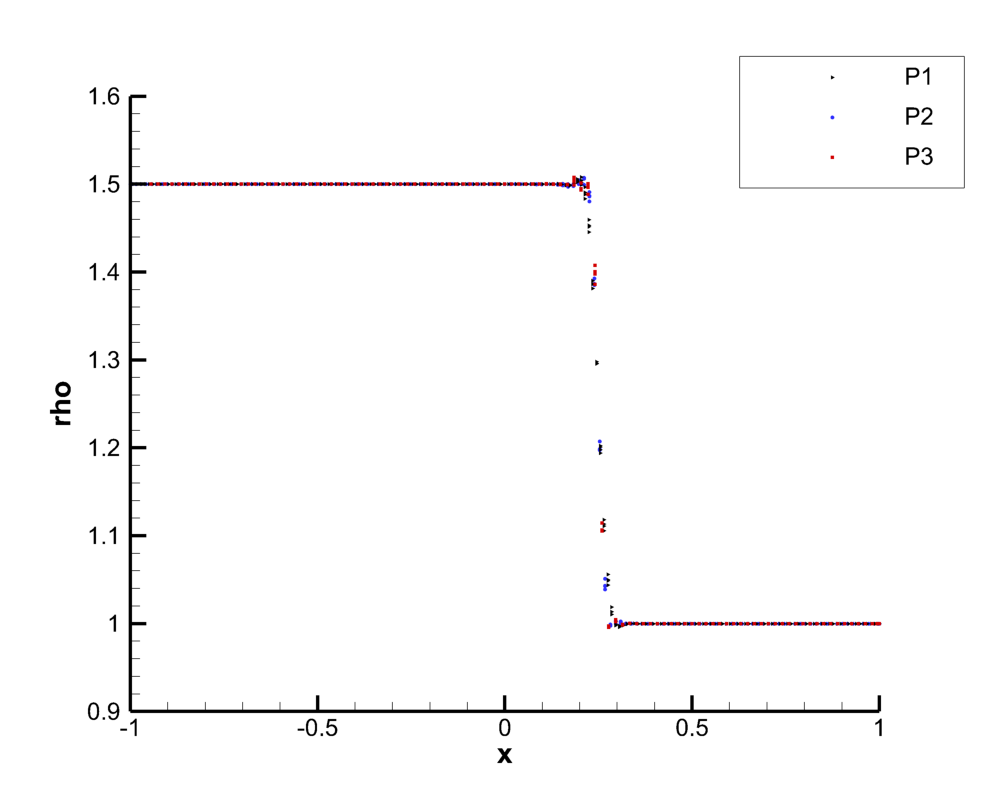}
	\includegraphics[width=0.49\linewidth]{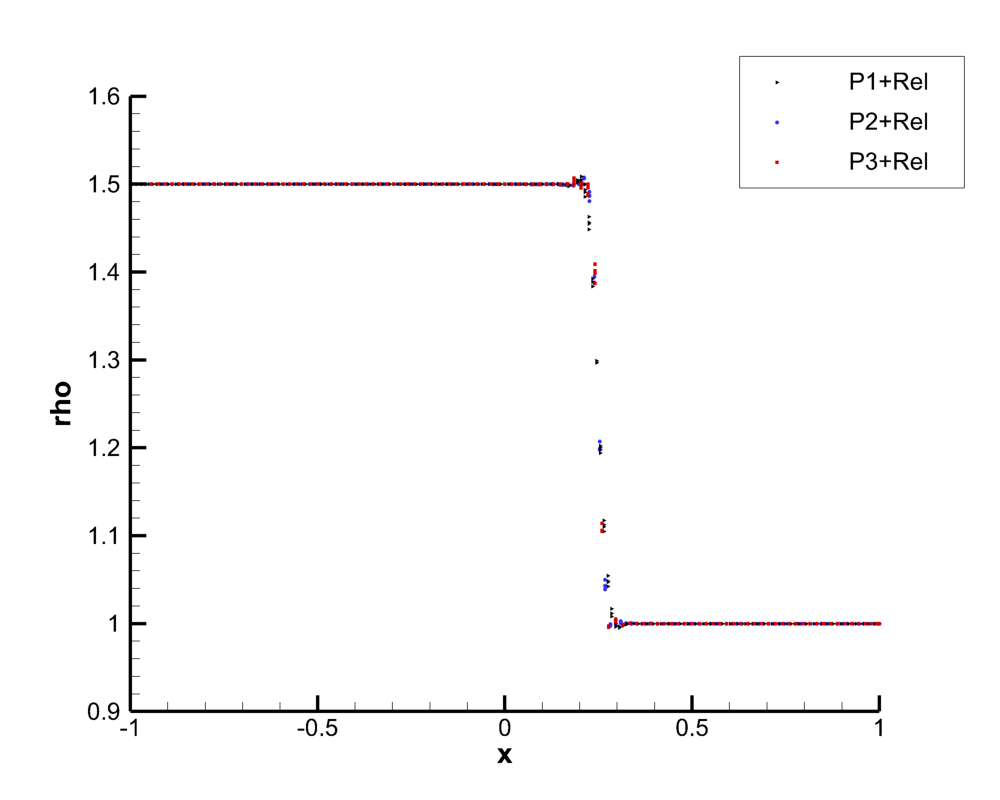}
	\caption{ Simulation of contact discontinuity traveling in the domain. On the left we used ADER-DG, on the right ADER-DG with entropy correction and relaxation method.
}
	\label{fig:eulercontactend}
\end{figure}

In \cref{fig:eulercontactend} we show the final solution for the ADER-DG (left) and the relaxation ADER-DG (right). We run the $\mathbb P_1$ scheme over a mesh of $5516$ elements ($16548$ DOFs), 
the $\mathbb P_2$ scheme over a mesh of $2860$ elements ($17160$ DOFs), 
and the $\mathbb P_3$ scheme over a mesh of $1744$ elements ($17440$ DOFs). We remark that both schemes do not have a limiter suited for shocks, hence small overshoots are visible at the sides of the contact discontinuity. 
The difference between the two simulations is not so evident, but, again, in classical methods there is a loss of total entropy.
\begin{figure}
	\centering
	\includegraphics[width=1.0\linewidth]{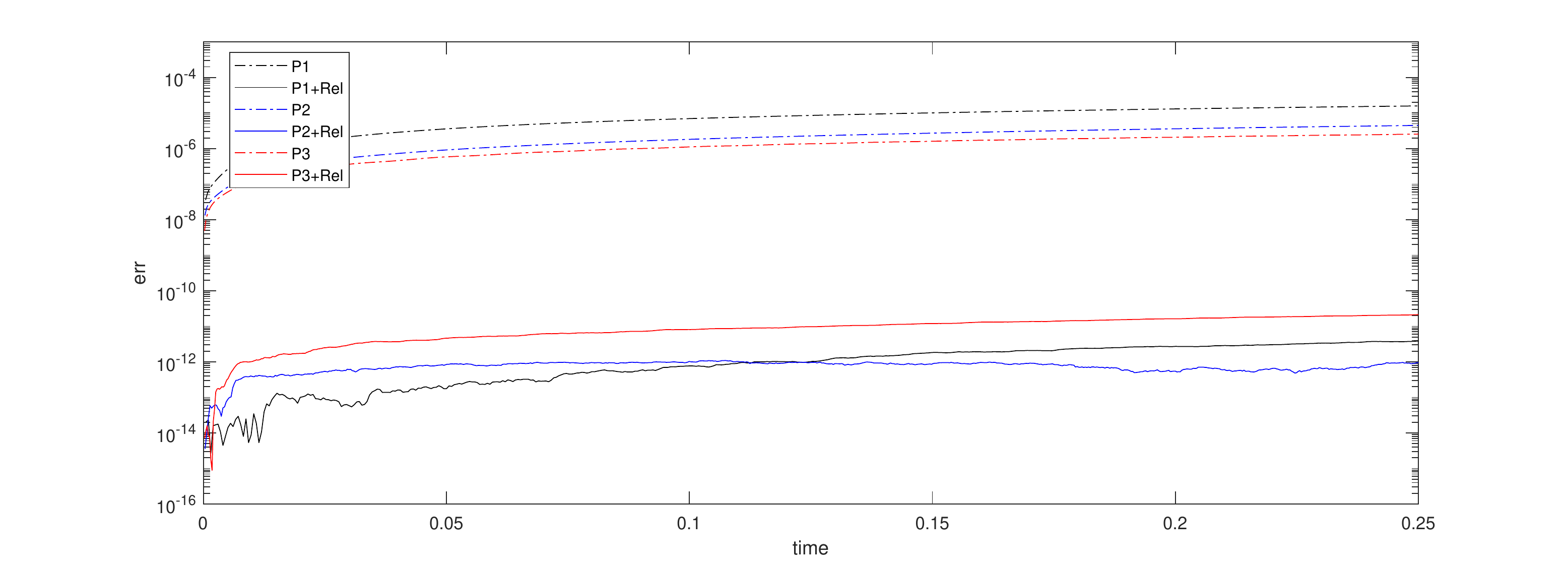}
	\caption{Entropy error for contact discontinuity test case.}
	\label{fig:eulercontactentropy}
\end{figure}
\begin{figure}
	\centering
	\includegraphics[width=0.5\linewidth]{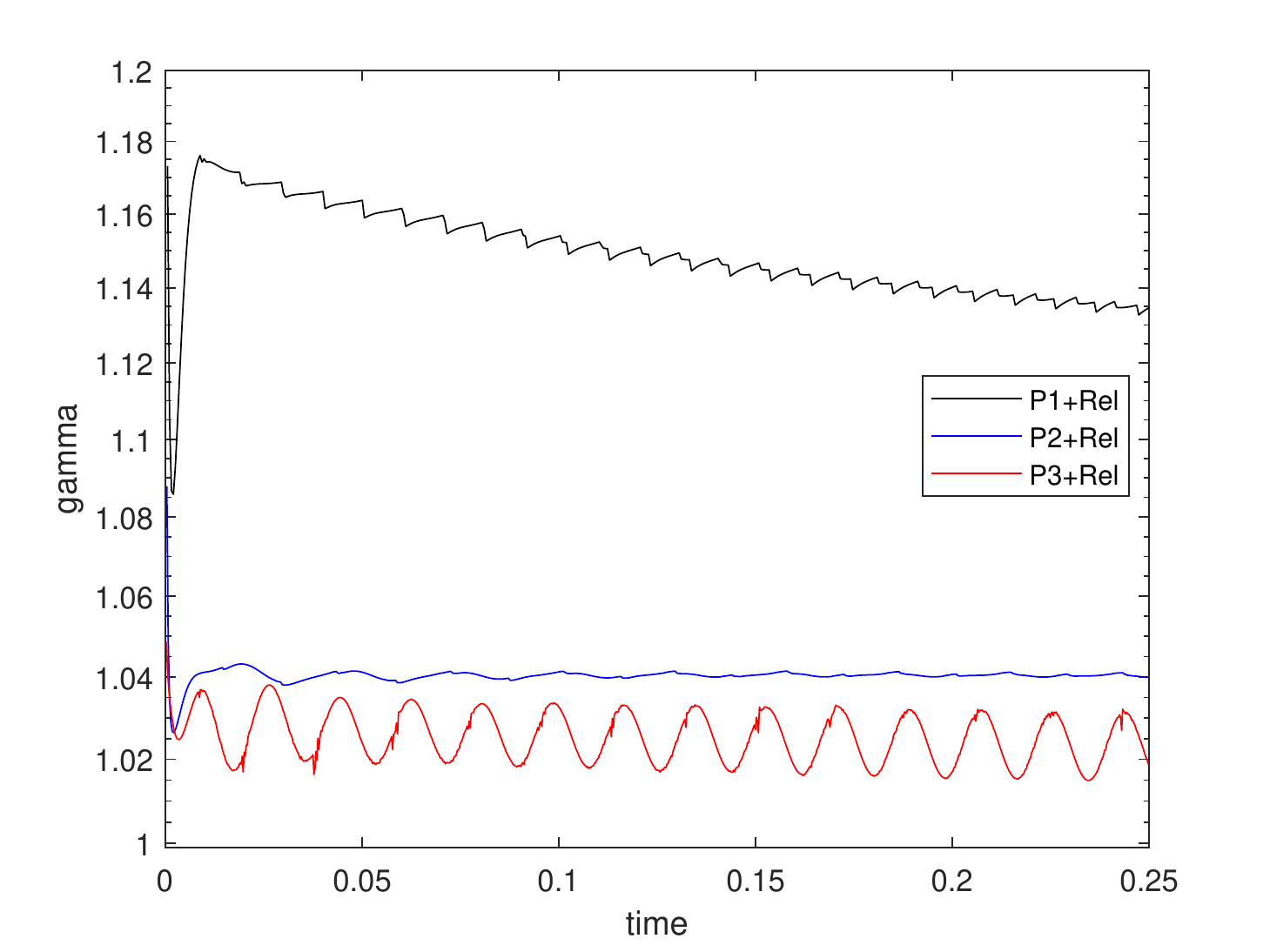}
	\caption{Values of $\gamma^n$ in time with relaxation algorithms for different orders but similar number of DOFs.}
	\label{fig:eulercontactgamma}
\end{figure}
In \cref{fig:eulercontactentropy} the total entropy of classical methods is lower than the exact one, and it is very similar for all the degrees of freedom. Here, we do not see huge advantages with a high order method as the solution is discontinuous. On the other side, for relaxation methods the exact value of the total entropy is preserved up to machine precision. 
Finally, in \cref{fig:eulercontactgamma} we show the values of $\gamma^n$ in time and we can see that they are much larger than in previous tests. This is again due to the discontinuous solution that makes the accuracy of the scheme decrease, needing a larger modification to obtain the total entropy conservation.

\subsubsection{123 Problem}
The last test that we propose is a two dimensional generalization of the 123 problem also presented in~\cite[Section 4.3.3]{toro-book}. It is defined on a domain $\Omega = [-L,L]^2$ with $L=1.2$, with initial condition
\begin{equation}
	\begin{pmatrix}
		\rho \\ u \\v \\ p
	\end{pmatrix} = 
\begin{pmatrix}
\rho_0 \\ u \\v \\ p_0
\end{pmatrix} = 
\begin{pmatrix}
1.0 \\ v_0 \frac{ \vec{x}}{||\vec{x}|| + 1e-4} \\ 0.4
\end{pmatrix},
\end{equation}
with $v_0=2.0$. We impose transmissive boundary conditions on all the domain.


\begin{figure}
	\centering
	\includegraphics[width=0.49\linewidth]{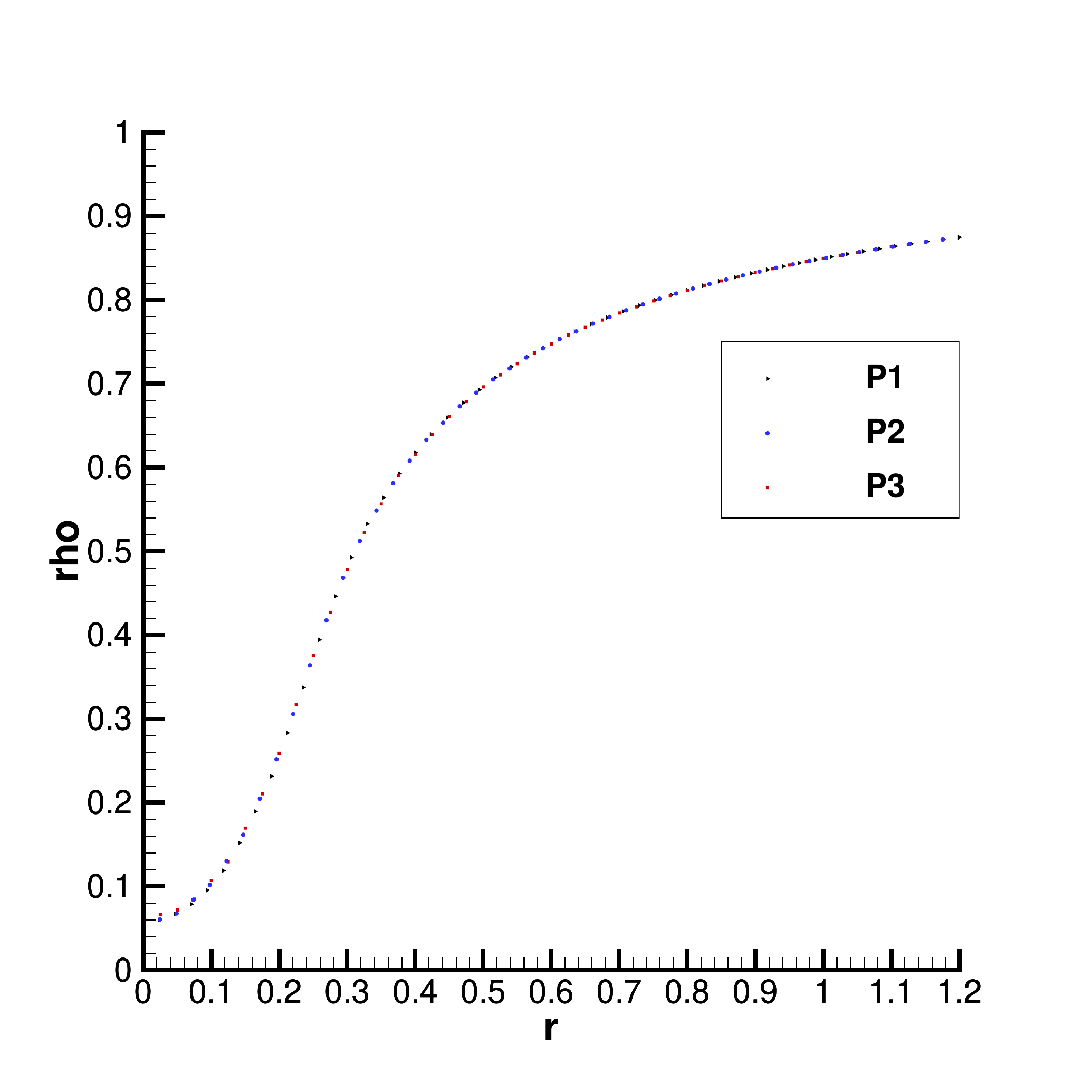}
	\includegraphics[width=0.49\linewidth]{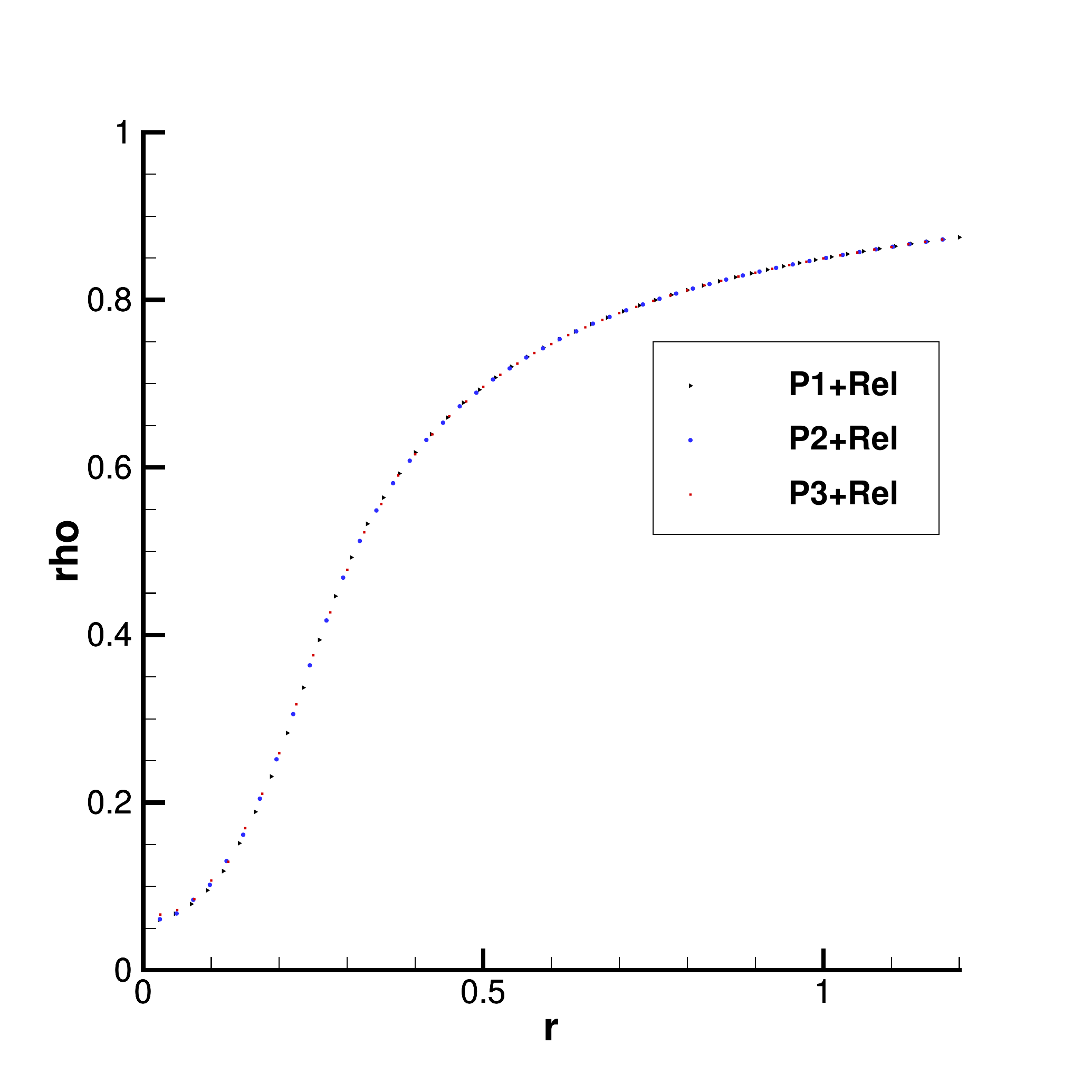}
	\caption{Simulation of the 123 problem. On the left we used ADER-DG, on the right ADER-DG with entropy correction and relaxation method.
	}
	\label{fig:eulerrarefactioncomparison}
\end{figure}

First of all, we plot the profile of the density solution at the final time $t_f=0.15$ on the line $y=0$. In \cref{fig:eulerrarefactioncomparison} we compare the classical ADER-DG and the relaxation ADER-DG methods. Again the different orders are run on different meshes so that the DOFs are more or less similar. We run the $P_1$ scheme over a mesh of $10108$ elements ($30324$ DOFs), the $P_2$ scheme over a mesh of $5130$ elements ($30780$ DOFs), and the $P_3$ scheme over a mesh of $3388$ elements ($33880$ DOFs). No visible differences are observable. 

\begin{figure}
	\centering
	\includegraphics[trim= 10 10 10 10,clip,width=0.49\linewidth]{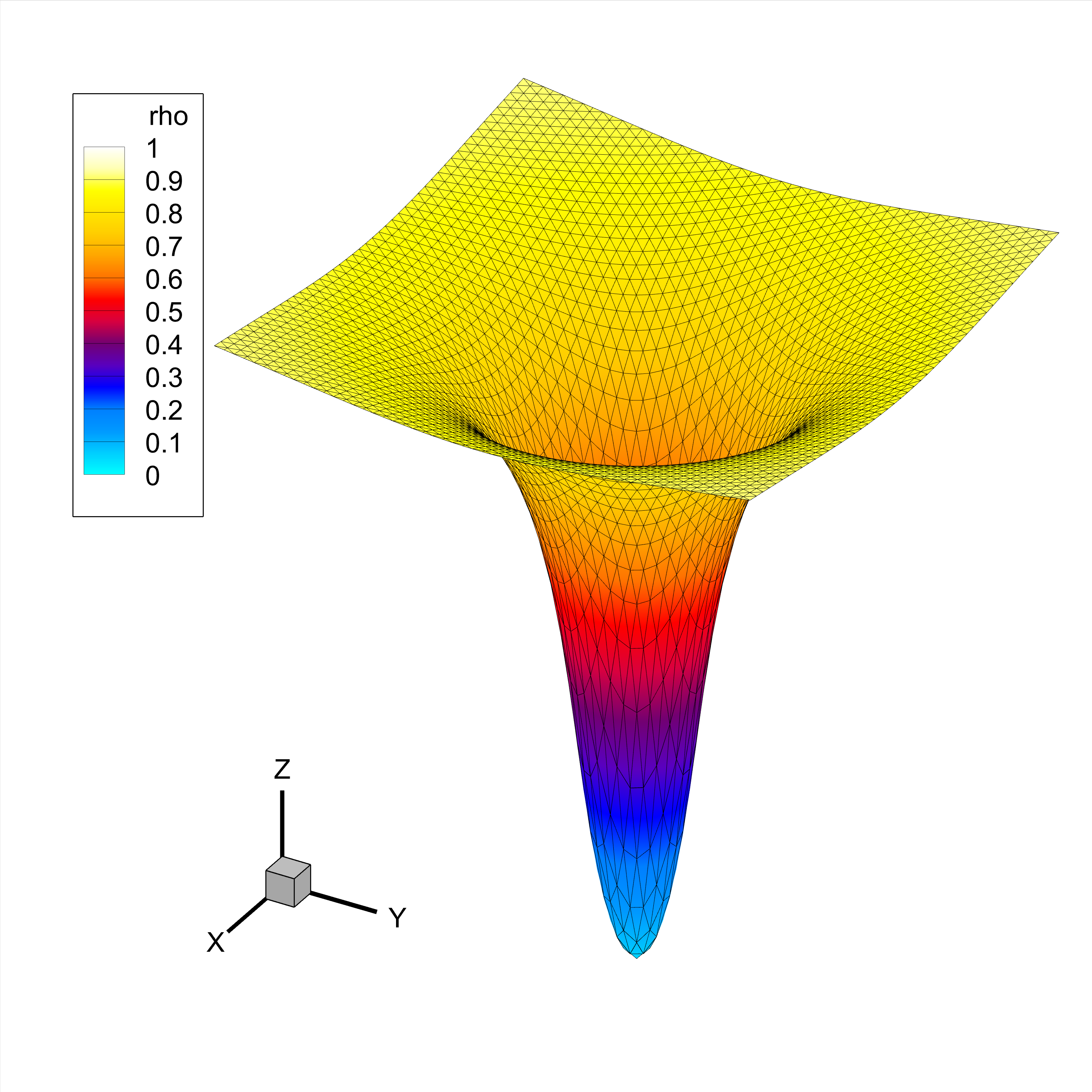}
	\includegraphics[trim= 10 10 10 10,clip,width=0.49\linewidth]{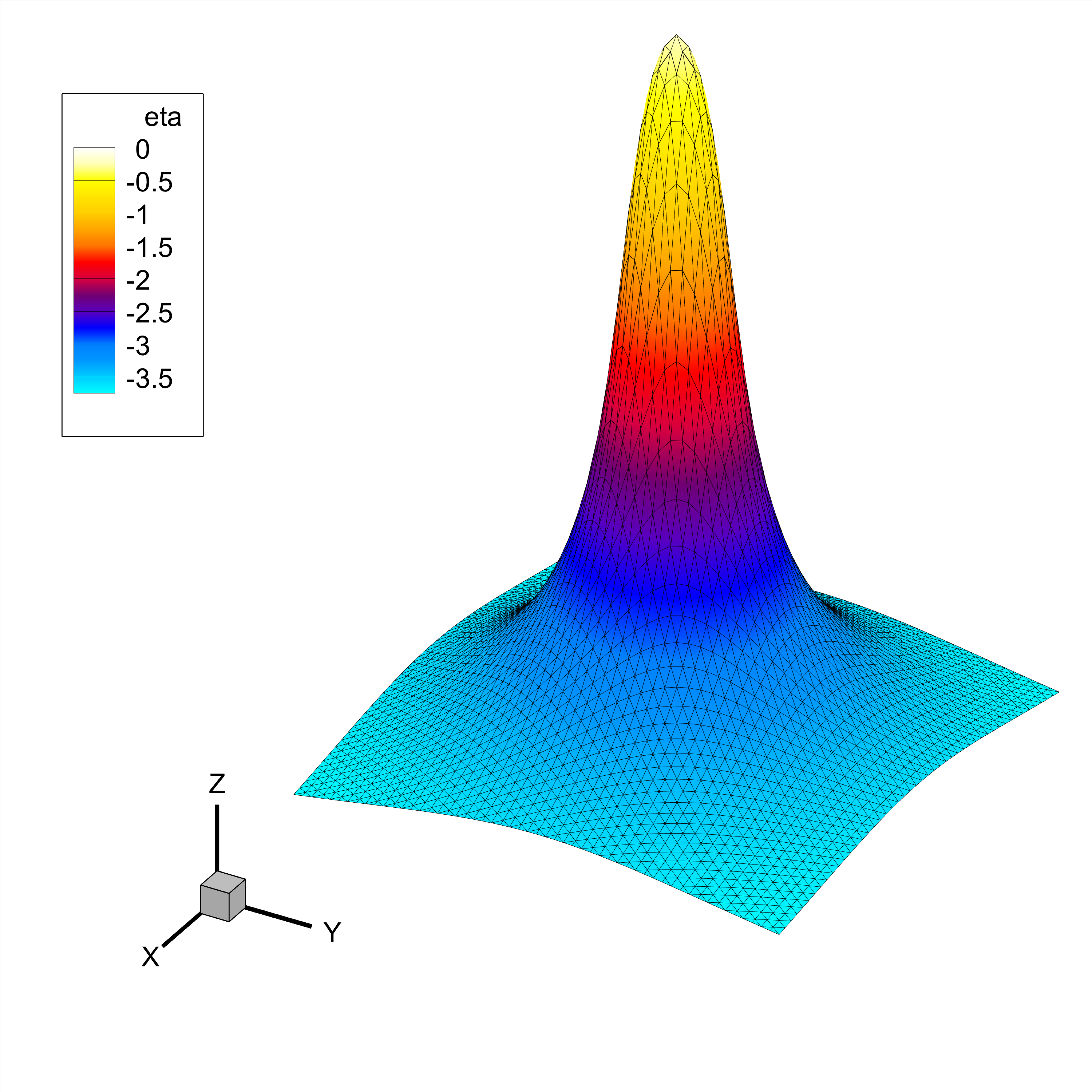}
	\caption{123 Problem.  
		In the figure we report the profile of the density $\rho$ (left) and the entropy $\eta$ (right) at time $t=0.15$ obtained with a DG scheme of order $3$ with relaxation on a very coarse mesh made of $5130$ elements.
	}
	\label{fig:eulerrarefactionplot}
\end{figure}
\begin{figure}
	\centering
	\includegraphics[width=1.0\linewidth]{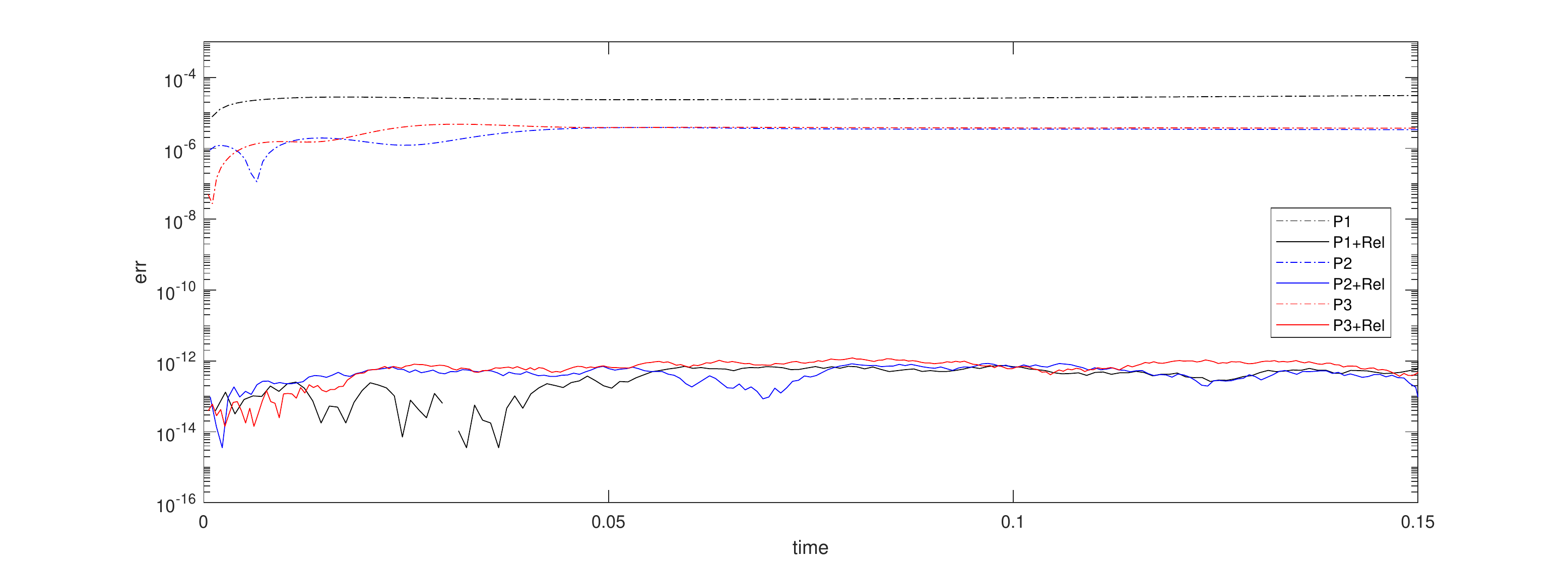}
	\caption{Entropy error for the 123 problem.}
	\label{fig:eulerrarefactionentropy}
\end{figure}
\begin{figure}
	\centering
	\includegraphics[width=0.5\linewidth]{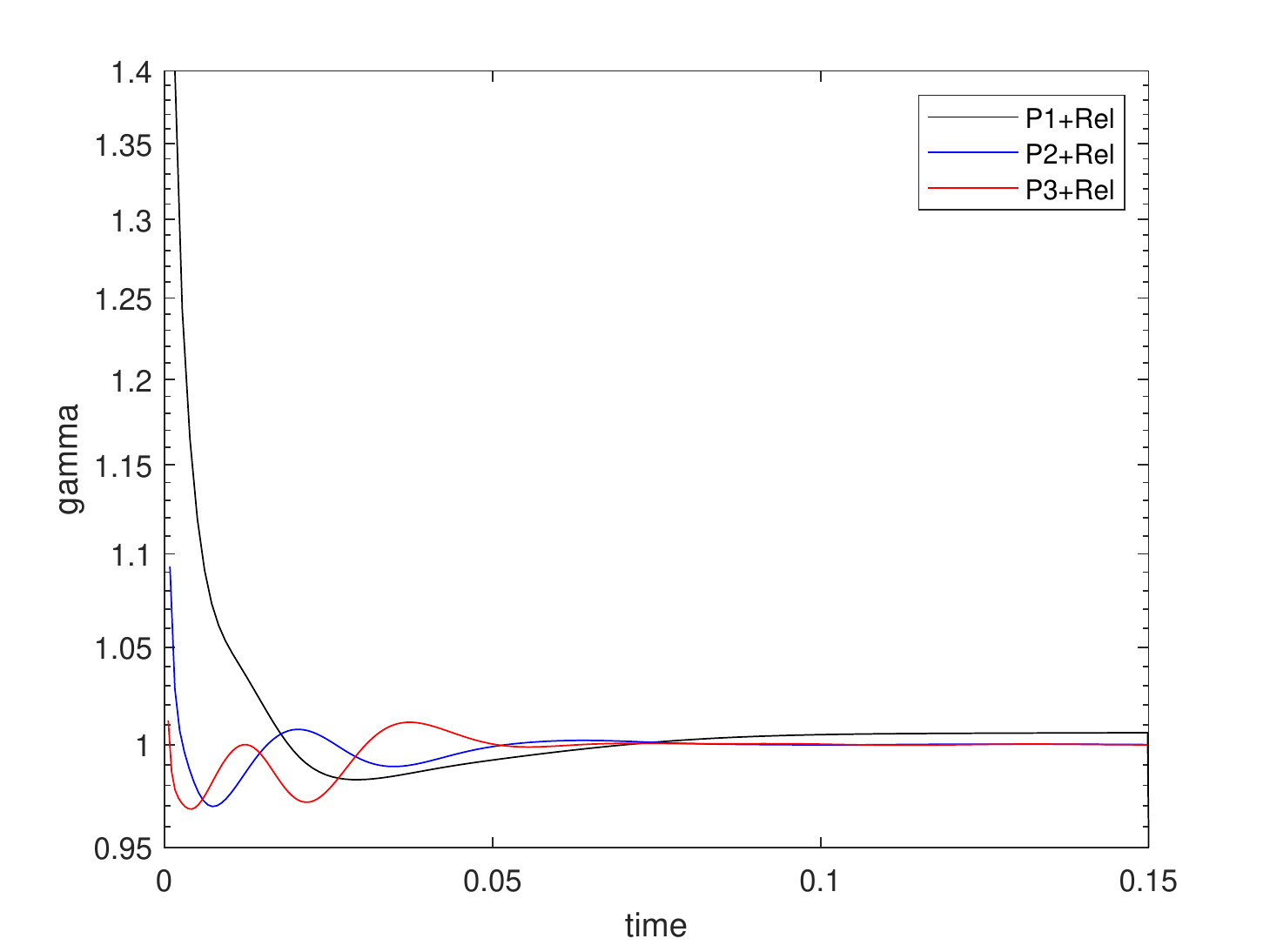}
	\caption{Values of $\gamma^n$ in time with relaxation algorithms for different orders but similar number of DOFs.}
	\label{fig:eulerrarefactiontgamma}
\end{figure}

In \cref{fig:eulerrarefactionplot} we plot again the solution $\rho$ and $\eta(\u)$ at the final time $t_f=0.15$ on the whole domain with a $P_3$ ADER-DG scheme on a coarse mesh with $5130$ elements. In this case, the boundary is affected by the propagation of the solution, hence, we do not know an analytical solution even for the total entropy.
Hence, to obtain the error of the total entropy error in \cref{fig:eulerrarefactionentropy} we compute numerically the integral of the entropy flux coming out at the boundary of the domain to approximate the exact value of the total entropy, i.e.,
\begin{equation}
	\mathcal{E}(\u(t^n))\approx \mathcal{E}(\u(t^{n-1})) - \sum_{s=0}^S \beta^s \dt^n \gamma^n \int_{\partial \Omega} \G(\u_h^{(s)}) \cdot \n\, \diff S.
\end{equation}
In \cref{fig:eulerrarefactionentropy} we observe again the usual behavior, the classical schemes preserve the total entropy only up to the accuracy of the method and, since the solution starts with a discontinuity in the velocity field, there is no large difference between $P_2$ and $P_3$. For the relaxation ADER-DG we are able to preserve the exact total entropy up to machine precision for all the orders.

Finally, in \cref{fig:eulerrarefactiontgamma} we plot the $\gamma^n$ values and we can observe that there is no large difference between $P_2$ and $P_3$, as seen already for the classical methods.

\section{Conclusions} \label{sec.concl}
In this work we have presented an enhancement for the ADER-DG arbitrarily high order method to provably conserve or dissipate the total entropy of the solution. 
The method is based on two ingredients. 
The first is an entropy correction term that allows to balance the spatial entropy production/destruction in each cell of the domain. Moreover, with the ADER-DG discretization is then possible to split the spatial discretization into an entropy conservative term and an entropy diffusive one.
The second ingredient is the relaxation method, which exploit the split spatial discretization and it is able to guarantee a zero entropy production also during the time discretization for the entropy conservative terms.
In practice, the proposed scheme outperforms the classical ADER-DG as it is able to exactly preserve the total entropy for many tests where physically this happens.
Indeed, we have tested the method on smooth solutions and a contact discontinuity and in all cases the method has proved its properties.

The method as it is has still some limitations as we started from an unlimited version of the ADER-DG. 
In order to be able to deal also with shocks and more challenging tests, it will be necessary to introduce a limiter and to insert it in the spatial discretization split so that also in the relaxation process the dissipation will occur only when needed.
Another extension of the method on which the authors are working on, is the $\mathbb P^N  \mathbb P ^M$ method~\cite{gaburro2021posteriori}, to obtain entropy conservative/dissipative schemes.

\section*{Acknowledgment}
E.~G. and M.~R. are members of the CARDAMOM team at the Inria center of the university of Bordeaux, D.~T. has been funded by a postdoctoral fellowship in Inria Bordeaux at the team CARDAMOM and by a postdoctoral fellowship in SISSA Trieste.
P.~O. gratefully acknowledge the support of the Gutenberg Research College and also wants to thank M.~R. for his invitation to Inria Bordeaux.  \\
E.~G. gratefully acknowledges the support received from the European Union’s Horizon 2020 Research and Innovation Programme under the Marie Sk\l{}odowska-Curie Individual Fellowship \textit{SuPerMan}, grant agreement No. 101025563.

\bibliographystyle{plain}
\bibliography{references}

\end{document}